\documentclass[leqno,12pt]{article}

\usepackage{latexsym}

\newtheorem{definition}{Definition}
\newtheorem{theorem}{Theorem}
\newtheorem{proposition}{Proposition}
\newtheorem{corollary}{Corollary}

\newcommand{\Prf}{Proof}

\makeatletter
   \renewcommand{\theequation}{%
   \thesection.\arabic{equation}}
   \@addtoreset{equation}{section}
\makeatother

\usepackage{graphicx,color}
\usepackage[T1]{fontenc}
\usepackage{textcomp}
\newcommand{\sq}{\textquotesingle}
\newcommand{\vc}[1]{\mbox{\boldmath{$#1$}}}

\newcommand{\mfill}[1]{\makebox[#1\textwidth][l]{}}
\newcommand{\LB}{\left\{}
\newcommand{\RB}{\right\}}
\newcommand{\LA}{\left(}
\newcommand{\RA}{\right)}

\newcommand{\TL}{\tilde}
\newcommand{\tvc}[1]{\TL{\mbox{\boldmath{$#1$}}}}
\newcommand{\tA}{\TL{A}}

\newcommand{\T}{{\rm T}}
\newcommand{\PO}{M}
\newcommand{\Pinv}{\PO^{-1}}
\newcommand{\Pinvt}{\PO^{-\T}}
\newcommand{\SH}{\sharp}
\newcommand{\NON}{\nonumber}
\newcommand{\UB}{\underline{ }}
\newcommand{\pR}{R}
\newcommand{\pP}{P}

\newcommand{\pZ}{\TL{\lambda}}
\newcommand{\PCONV}{=}
\def\L({\left(}
\def\R){\right)}

\def\km{{k-1}}
\def\P{M}
\def\KSP2{{\mathcal K}_{k+1}}
\def\KSP{${\mathcal K}_{k+1}$}

\def\NON{\nonumber}
\def\SP{\;}  % thick space
\def\UL{\underline}
\def\UB{\underline{ }}
\newcounter{alg}
\newenvironment{indention}[1]%
{\par\begingroup\addtolength{\leftskip}{#1}}%
{\par\endgroup}
\def\UL{\underline}
\def\UB{\underline{ }}
\newcommand{\ALGTITLEWIDTH}{2em}

\newcommand{\ALGWIDTHB}{6em}

\def\FIGSIZE{0.7}
\def\ALGCHARTSIZE{0.8}
\newcommand{\Def}{Definition}
\newcommand{\Dfn}{\Def}

\newcommand{\Thm}{Theorem}

\newcommand{\Prop}{Proposition}
\newcommand{\prop}{proposition}
\newcommand{\Cor}{Corollary}
\newcommand{\cor}{corollary}

\newcommand{\Alg}{Algorithm} 
\newcommand{\Ref}{Reference} 
\newcommand{\Sec}{Section}

\newcommand{\App}{Appendix}
\newcommand{\Fig}{Figure}
\newcommand{\Tab}{Table }
\newcommand{\reeqno}[1]{\renewcommand{\theequation}{#1}}
\def\PCGS{PCGS}
\def\GRAPHARR{.}
\def\GRAPHTRR{.}
\def\GRAPHTRE{.}
\newcommand{\GRAPHisrv}{.}
\def\GRAPHDIRTWO{.}
\title{
Structure of the preconditioned system
in various preconditioned conjugate gradient squared algorithms
}

\author{
Shoji Itoh%
\thanks{Department of Engineering Science,
     Faculty of Engineering,
     Osaka Electro-Communication University. 
(itosho@acm.org).
}
$\;\;$
and
$\;$
Masaaki Sugihara%
\thanks{Department of Physics and Mathematics,
     College of Science and Engineering, Aoyama Gakuin University.
}
\thanks{Deceased 5 January 2019}
}

\date{}

\begin{document}
\maketitle

%%%%% Begin Abstract %%%%%%%%%%%
\begin{abstract}
An improved preconditioned conjugate gradient squared (PCGS) algorithm
has recently been proposed, and
it performs much better than the conventional PCGS algorithm.
In this paper,
the improved PCGS algorithm is verified as a coordinative to
the left-preconditioned system,
and it has the advantages of both the conventional and the left-PCGS;
this is done
by comparing, analyzing,
and executing numerical examinations of various PCGS algorithms,
including another improved one.
We show that
the direction of the preconditioned system for the CGS method is determined
by the operations of $\alpha_k$ and $\beta_k$ in the PCGS algorithm.
By comparing the logical structures of these algorithms,
we show that the direction of the preconditioned system can be switched
by the construction and setting of the initial shadow residual vector.
\\

\noindent
{\bf Keywords:}
Congruence of preconditioning conversion,
Direction of preconditioned system,
Preconditioned Krylov subspace,
Improved PCGS
\end{abstract}
%%%%% end %%%%%%%%%%%

%%%% Start %%%%%%
\section{Introduction}
\label{sec:intro}

The conjugate gradient squared (CGS) \cite{sonneveld1989}
is one of various methods used to solve systems of linear equations
\begin{eqnarray}
A\vc{x} &=& \vc{b},
\label{eqn:linear}
\end{eqnarray}
where
the coefficient matrix $A$ of size $n\times n$ is usually nonsymmetric,
$\vc{x}$ is the solution vector,
and
$\vc{b}$ is the right-hand side (RHS) vector.

The CGS is a bi-Lanczos method
that belongs to the class of Krylov subspace methods.
Bi-Lanczos-type methods are derived from
the bi-conjugate gradient (BiCG) method \cite{fletcher1976,lanczos1952},
which assumes the existence of a dual system
$A^\T\vc{x}^\SH = \vc{b}^\SH$ %.
(we will refer to this as the ``shadow system'').
Bi-Lanczos-type algorithms have the advantage of requiring less memory than do
Arnoldi-type algorithms, which is another class of Krylov subspace methods.
Furthermore, a variety of bi-Lanczos-type algorithms,
such as the bi-conjugate gradient stabilized (BiCGStab) method \cite{vorst1992}
and the generalized product-type method based on the BiCG (GPBiCG) \cite{zhang1997},
have been constructed by adopting the idea behind the derivation of the CGS.
Various iterative methods, including bi-Lanczos-type algorithms,
are often used following a preconditioning operation
that is used to improve the properties of the linear equations.
Such algorithms are called preconditioned algorithms;
for example, the preconditioned CGS (PCGS).
Therefore, it is very important to study the properties of the PCGS so that its performance can be improved.

Generally,
the degree $k$ of the Krylov subspace generated
by $A$ and $\vc{r}_0$ is expressed as
${\cal K}_k\L(A,\vc{r}_0\R)$
$={\rm span}\LB\vc{r}_0, A\vc{r}_0, A^2\vc{r}_0, \cdots, A^{k-1}\vc{r}_0\RB$,
where
$\vc{r}_0$ is the initial residual vector $\vc{r}_0 = \vc{b}-A\vc{x}_0$, and $\vc{x}_0$
is the initial guess at the solution.
The Krylov subspace ${\cal K}_k\L(A,\vc{r}_0\R)$ generated by the $k$-th iteration
forms the structure of
$\vc{x}_k \in \vc{x}_0 + {\cal K}_k\L(A,\vc{r}_0\R)$, %\NON
where $\vc{x}_k$ is the approximate solution vector
(or simply the ``solution vector'').
However,
for a given preconditioned Krylov subspace method,
there are various different algorithms that can be used for the preconditioning conversion.
In such cases,
the structure of the approximate solution formed by
the Krylov subspace is often different for different algorithms,
and
the performance of these various algorithms can also differ substantially \cite{itoh2015a}.

An improved PCGS algorithm has been proposed \cite{itoh2015a}.
\Ref~\cite{itoh2015a} shows that this improved algorithm has
many advantages over the conventional PCGS algorithms \cite{barrett1994, meurant2005, vorst1992}.
In this paper,
a variety of PCGS algorithms are discussed.
We begin by considering
two typical PCGS algorithms, and
we analyze the mathematical structure
with respect to and in connection with the Krylov subspace.
In particular,
we define and consider the direction of the preconditioned system of these PCGS algorithms.
We then perform the same analysis for two improved PCGS algorithms,
one of which was mentioned above \cite{itoh2015a},
and the other is presented in the present paper.
However,
we note that
it is not our purpose to propose a new algorithm,
but to analyze the improved PCGS algorithms
by comparing the logical structures
and the numerical results
of four different PCGS algorithms.

In this paper, when we refer to a
{\it preconditioned algorithm}, we mean one that uses a
preconditioning operator $\PO$ or a preconditioning matrix,
and by {\it preconditioned system}, we mean
one that has been converted by some operator(s) based on $\PO$.
These terms never indicate the
{\it algorithm for the preconditioning operation itself},
such as incomplete LU decomposition or by using the approximate inverse.
For example,
under a preconditioned system,
the original linear system (\ref{eqn:linear}) becomes
\begin{eqnarray}
&&
\tA\tvc{x} = \tvc{b},
\label{eqn:plinear}
\\
&&
\tA     \PCONV \Pinv_L A \Pinv_R, \;\;
\tvc{x} \PCONV \PO_R\vc{x}, \;\;
\tvc{b} \PCONV \Pinv_L\vc{b},
\label{eqn:plinear_conv}
\end{eqnarray}
with the preconditioner $\PO = \PO_L\PO_R$ ($\PO \approx A$).
In this paper,
the matrix and the vector under the preconditioned system
are indicated by a tilde ($\;\TL{ }\;$).
However,
the conversions in (\ref{eqn:plinear}) and (\ref{eqn:plinear_conv}) are
not implemented directly;
rather, we construct the preconditioned algorithm
that is equivalent to solving (\ref{eqn:plinear}).

This paper is organized as follows.
\Sec~\ref{sec:pcgs_alg} provides
various preconditioned CGS algorithms.
In particular, we consider
the right- and left-directions of the preconditioned systems for CGS algorithms.
The improved PCGS algorithms are shown to be
coordinative to the left-preconditioned system.
\Sec~\ref{sec:congruent_and_direction_pcgs} discusses the difference
between
the direction of a preconditioning conversion 
and
the direction of a preconditioned system.
We show that preconditioning conversions are congruent for PCGS algorithms, 
and
we provide some examples in which the direction of the preconditioned system for the CGS is switched.
In Section~\ref{sec:numerical_experiments},
we present some numerical results to illustrate
the convergence properties
of the various PCGS algorithms discussed in Section~\ref{sec:pcgs_alg},
and we illustrate the effect of
switching the direction of the preconditioned system for the CGS algorithm 
in Section~\ref{sec:congruent_and_direction_pcgs}.
Finally,
our conclusions are presented in Section~\ref{sec:conclusion}.

\section{Analyses of various PCGS algorithms}
\label{sec:pcgs_alg}

In this section, four kinds of PCGS algorithms are analyzed.
These PCGS algorithms can be derived as follows.
\\

%%%%%%%%%%%%%%%%%%%%%%%%%%%%%%%%%%%%%%%%%%%%%%%%%%%%%%%%%%%%%
\refstepcounter{alg}
\label{alg:pcgs_simple}
\noindent
\hspace{\ALGTITLEWIDTH}
{\bf Algorithm \thealg.
CGS under preconditioned system:}
\begin{indention}{\ALGWIDTHB}
%%%%%%%%%%%%%%%%%%%%%%%%%%%%%%%%%%%%%%%%%%%%%%%%%%%%%%%%%%%%%
\noindent
$\tvc{x}_0$ is the initial guess,
$\SP\tvc{r}_0= \tvc{b}-\tA\tvc{x}_0$,
set $\SP\beta^{\rm \PCGS}_{-1}=0$,
\\
$\LA \tvc{r}^\SH_0, \tvc{r}_0 \RA \neq 0$,
e.g.,
$\tvc{r}^\SH_0 = \tvc{r}_0, \SP $
 \\
For $k = 0, 1, 2, \cdots ,$ until convergence, Do:
\begin{eqnarray}
&& \tvc{u}_k = \tvc{r}_k + \beta^{\rm \PCGS}_{k-1}\tvc{q}_{k-1}, \NON \\
&& \tvc{p}_k = \tvc{u}_k + \beta^{\rm \PCGS}_{k-1}\L(\tvc{q}_{k-1} + \beta^{\rm \PCGS}_{k-1}\tvc{p}_{k-1}\R), \NON
\\
&& \alpha^{\rm \PCGS}_k = \frac
 {\LA \tvc{r}^\SH_0,    \tvc{r}_k \RA}
 {\LA \tvc{r}^\SH_0, \tA\tvc{p}_k \RA}, \NON
\\
&& \tvc{q}_k     = \tvc{u}_k - \alpha^{\rm \PCGS}_k\tA   \tvc{p}_k, \NON \\
&& \tvc{x}_{k+1} = \tvc{x}_k + \alpha^{\rm \PCGS}_k   \L(\tvc{u}_k + \tvc{q}_k\R), %\NON \\
\label{eqn:x_pcgs}
\\
&& \tvc{r}_{k+1} = \tvc{r}_k - \alpha^{\rm \PCGS}_k\tA\L(\tvc{u}_k + \tvc{q}_k\R), \NON
\\
&& \beta^{\rm \PCGS}_k = \frac
 {\LA \tvc{r}^\SH_0, \tvc{r}_{k+1} \RA}
 {\LA \tvc{r}^\SH_0, \tvc{r}_k     \RA}, \NON
\end{eqnarray}
End Do
\\

\end{indention}
%%%%%%%%%%%%%%%%%%%%%%%%%%%%%%%%%%%%%%%%%%%%%%%%%%%%%%%%%%%%%

Any preconditioned algorithm can be derived by substituting
the matrix with the preconditioner for the matrix with the tilde
and
the vectors with the preconditioner for the vectors with the tilde.
Obviously,
\Alg~\ref{alg:pcgs_simple} without the preconditioning conversion
is the same as the CGS.
If $\tA$ is a symmetric matrix and
$\tvc{r}^\SH_0=\tvc{r}_0$,
then \Alg~\ref{alg:pcgs_simple} can be adapted while maintaining its symmetric property.

The case shown in (\ref{eqn:plinear_conv}) is called
{\it two-sided preconditioning},
the case in which $\PO_L = \PO$ and $\PO_R = I$ is called {\it left preconditioning},
and the case in which $\PO_L = I$ and $\PO_R = \PO$ is called {\it right preconditioning},
where $I$ denotes the identity matrix.
We now formally define these%
\footnote{
Here,
we have offered a general definition.
However,
for preconditioned bi-Lanczos-type algorithms,
additional restrictions are necessary \cite{itoh2016b}.
}%
.

\begin{definition}
\label{dfn:dir_prec}
For the system and solution
\vspace*{10pt}
\begin{equation}
\reeqno{\ref{eqn:plinear}\sq}
% \hspace*{-148pt}
\hspace*{-167pt}
\tA\tvc{x} = \tvc{b},
\end{equation}
% \vspace*{-30pt}
\begin{equation}
\reeqno{\ref{eqn:plinear_conv}\sq}
\tA     \PCONV \Pinv_L A \Pinv_R, \;\;
\tvc{x} \PCONV \PO_R\vc{x}, \;\;
\tvc{b} \PCONV \Pinv_L\vc{b},
\end{equation}
we define
the direction of a preconditioned system of linear equations as follows:
\begin{itemize}
\item {\it The two-sided preconditioned system:}
Equation (\ref{eqn:plinear_conv}\sq);
\vspace{-5pt}
\item {\it The right-preconditioned system:}
$\PO_L = I$ and $\PO_R = \PO$ in  (\ref{eqn:plinear_conv}\sq);
\vspace{-5pt}
\item {\it The left-preconditioned system:}
$\PO_L = \PO$ and $\PO_R = I$ in  (\ref{eqn:plinear_conv}\sq),
\end{itemize}
where
$\PO$ is the preconditioner $\PO=\PO_L\PO_R$ ($\PO \approx A$), and
$I$ is the identity matrix.

Other vectors in the solution method are not preconditioned.
The initial guess is given as $\vc{x}_0$,
and $\tvc{x}_0 = \PO_R \vc{x}_0$.
\end{definition}

\setcounter{equation}{1}

The two-sided preconditioned system may be impracticable, but it is of theoretical interest.

The preconditioned system is different from the preconditioning conversion
as the next definition.
There are various ways of performing a preconditioning conversion,
but the direction of the preconditioned system is uniquely defined.

\begin{definition}[Congruence]
\label{dfn:congruence}

Let the term ``all directions on preconditioning conversions'' be
a synthesis of the preconditioning conversion for
``the two-sided direction'', ``the right direction'', and ``the left direction'',
not only for the system and solution
(\ref{eqn:plinear}\sq),
but also for other vectors in the solution method.

If all directions on the preconditioning conversions to the solution method are reduced to
one and the same algorithm description,
then we refer to this as ``congruence'' in the direction of the preconditioning conversion.
Furthermore,
the term ``congruency'' refers to the congruence property.
\end{definition}

As an example of congruence,
see the preconditioning conversion by
(\ref{eqn:conv_precb}),  (\ref{eqn:conv_precb2}), and (\ref{eqn:conv_precb1})
for Algorithm~\ref{alg:cgs_precb} in Section~\ref{sssec:conventional_pcgs}.

Both the CGS and the PCGS extend the two-dimensional subspace
in each iteration \cite{bruaset1995,gutknecht2000};
therefore, the
Krylov subspace ${\cal K}_{2k}(\tA,\tvc{r}_0)$ generated
by the $k$-th iteration forms the structure of
\begin{eqnarray}
\tvc{x}_k \in \tvc{x}_0 + {\cal K}_{2k}(\tA,\tvc{r}_0).
\label{eqn:ksp_cgs}
\end{eqnarray}

The CGS method is derived from the BiCG method \cite{sonneveld1989}.
Now,
the recurrence relations of the BiCG under the preconditioned system are
\begin{eqnarray}
\pR_0(\pZ) &=& 1, \SP\SP \pP_0(\pZ)=1, \NON
\label{eqn:prec_init_pol}
\\
\pR_k(\pZ) &=&
  \pR_{k-1}(\pZ) -\alpha^{\rm PBiCG}_{k-1}\pZ \pP_{k-1}(\pZ), 
\label{eqn:prec_res_pol}
\\
\pP_k(\pZ) &=&
  \pR_k(\pZ) +\beta^{\rm PBiCG}_{k-1} \pP_{k-1}(\pZ). 
\label{eqn:prec_prob_pol}
\end{eqnarray}
Here,
$\pR_k(\pZ)$ is the degree $k$ of the residual polynomial,
and
$\pP_k(\pZ)$ is the degree $k$ of the probing direction polynomial,
that is,
$\tvc{r}^{\rm BiCG}_k=\pR_k(\tA)\tvc{r}_0$
and
$\tvc{p}^{\rm BiCG}_k=\pP_k(\tA)\tvc{r}_0$.
The CGS method under preconditioned system can be derived
by introducing the idea of
$\{\pR_k(\pZ)\}^2$ and $\{\pP_k(\pZ)\}^2$ \cite{itoh2015a,sonneveld1989}.
Then,
we can represent the following relations between
these polynomials
and
the vectors in Algorithm \ref{alg:pcgs_simple}:
\begin{eqnarray}
&& \tvc{r}_k = \{\pR_k(\tA)\}^2\tvc{r}_0, \;\;\;\; \;\;\;\;
   \tvc{p}_k = \{\pP_k(\tA)\}^2\tvc{r}_0, %\;
\label{eqn:pol_pcgs}
\\
&& \tvc{u}_k =  \pP_k(\tA)\pR_k(\tA)\tvc{r}_0, \;\;\;\;
   \tvc{q}_k =  \pP_k(\tA)\pR_{k+1}(\tA)\tvc{r}_0.
\NON
\end{eqnarray}

\subsection{Two typical PCGS algorithms}
\label{ssec:typical_existing_pcgs}

In this subsection,
we present two well-known and typical PCGS algorithms.
One is a right-preconditioned system,
although this is not always recognized, and the other is a left-preconditioned system.
For each of these algorithms,
we examine the mathematical structure 
with respect to and in connection with the Krylov subspace and the solution vector.

\subsubsection{Conventional PCGS: Right-preconditioned CGS}
\label{sssec:conventional_pcgs}

This PCGS algorithm has been described in many manuscripts and numerical libraries;
for example,
 see \cite{barrett1994, meurant2005, vorst1992}.
It is usually derived by the following preconditioning conversion%
\footnote{%
In this case,
the initial shadow residual vector (ISRV) $\tvc{r}^\SH_0$
is converted to $\PO^\T_L\vc{r}^\flat_0$.
Here, 
the internal structure of $\vc{r}^\flat_0$ is
 $\vc{r}^\flat_0 \equiv \Pinvt\vc{r}^\SH_0$.
The notation $\vc{r}^\flat_0$ will be discussed
in Section~\ref{sec:congruent_and_direction_pcgs}.
The same applies to (\ref{eqn:conv_precb2}) and (\ref{eqn:conv_precb1}).
}%
:
\begin{eqnarray}
&&
  \tA        \PCONV \Pinv_L A \Pinv_R,  \SP\SP
  \tvc{x}_k  \PCONV \PO_R     \vc{x}_k, \SP\SP
  \tvc{b}    \PCONV \Pinv_L   \vc{b},   \SP\SP
  \tvc{r}_k     \PCONV \Pinv_L \vc{r}_k,       \SP\SP
\label{eqn:conv_precb}
\\
&&
  \tvc{r}^\SH_0 \PCONV \PO^\T_L  \vc{r}^\flat_0, \SP\SP
  \tvc{p}_k \PCONV \Pinv_L \vc{p}_k, \SP\SP
  \tvc{u}_k \PCONV \Pinv_L \vc{u}_k, \SP\SP
  \tvc{q}_k \PCONV \Pinv_L \vc{q}_k . %, \SP\SP
 \NON
\end{eqnarray}

Finally, Algorithm \ref{alg:cgs_precb} is derived.
% \\

\newpage
%%%%%%%%%%%%%%%%%%%%%%%%%%%%%%%%%%%%%%%%%%%%%%%%%%%%%%%%%%%%%
\refstepcounter{alg}
\label{alg:cgs_precb}
\noindent
\hspace{\ALGTITLEWIDTH}
{\bf Algorithm \thealg.
Conventional PCGS algorithm:}
\begin{indention}{\ALGWIDTHB}
%%%%%%%%%%%%%%%%%%%%%%%%%%%%%%%%%%%%%%%%%%%%%%%%%%%%%%%%%%%%%
\noindent
$\vc{x}_0$ is the initial guess,
$\SP\vc{r}_0= \vc{b}-A\vc{x}_0,\SP\SP$
set $\SP\beta_{-1}=0,$
\\
$\LA \TL{\vc{r}}^\SH_0, \TL{\vc{r}}_0\RA
 \PCONV \LA \vc{r}^\flat_0, \vc{r}_0\RA
 \neq 0$,
e.g.,
$\vc{r}^\flat_0 = \vc{r}_0,\SP$

\noindent
For $k = 0, 1, 2, \cdots ,$ until convergence, Do:
\begin{eqnarray}
 && \vc{u}_k = \vc{r}_k + \beta_{k-1}\vc{q}_{k-1}, \NON \\
 && \vc{p}_k = \vc{u}_k + \beta_{k-1}\L(\vc{q}_{k-1} + \beta_{k-1}\vc{p}_{k-1}\R), \NON \\
 && \alpha_k = \frac
             {\LA\vc{r}^\flat_0,  \vc{r}_k\RA}
             {\LA\vc{r}^\flat_0, A\Pinv\vc{p}_k\RA}, \NON
\label{eqn:conv_cgs_alpha}
\\
 && \vc{q}_k = \vc{u}_k - \alpha_k A\Pinv\vc{p}_k, \NON
\\
 && \vc{x}_{k+1} = \vc{x}_k + \alpha_k  \Pinv \L(\vc{u}_k + \vc{q}_k\R), %\NON \\
\label{eqn:sol_conv_pcgs}
\\
 && \vc{r}_{k+1} = \vc{r}_k - \alpha_k A\Pinv \L(\vc{u}_k + \vc{q}_k\R), \NON \\
 && \beta_k = \frac
             {\LA\vc{r}^\flat_0, \vc{r}_{k+1}\RA}
             {\LA\vc{r}^\flat_0, \vc{r}_k\RA}, \NON
\label{eqn:conv_cgs_beta}
\end{eqnarray}
End Do
\\

\end{indention}
%%%%%%%%%%%%%%%%%%%%%%%%%%%%%%%%%%%%%%%%%%%%%%%%%%%%%%%%%%%%%

The stopping criterion is
\begin{eqnarray}
&&\frac{\|\vc{r}_{k+1}\|}{\|\vc{b}\|} \leq \varepsilon .
\label{eqn:conv_judge}
\end{eqnarray}

The results of this algorithm can also be derived by the following conversion:
\begin{eqnarray}
&&
  \tA       \PCONV A\Pinv,\SP\SP
  \tvc{x}_k \PCONV  \PO\vc{x}_k,\SP\SP
  \tvc{b}   \PCONV  \vc{b},\SP\SP
  \tvc{r}_k     \PCONV \vc{r}_k,\SP\SP
\label{eqn:conv_precb2}
\\
&&
  \tvc{r}^\SH_0 \PCONV \vc{r}^\flat_0,\SP\SP %\NON
  \tvc{p}_k \PCONV \vc{p}_k,\SP\SP
  \tvc{u}_k \PCONV \vc{u}_k,\SP\SP
  \tvc{q}_k \PCONV \vc{q}_k . %,\SP\SP
 \NON
\end{eqnarray}

This is the same as using $\PO_L=I$ and $\PO_R=\PO $ in (\ref{eqn:conv_precb}).
Furthermore,
this is the same as converting
only $\tA$, $\tvc{x}_k$, and $\tvc{b}$,
that is, the right-preconditioned system.

Furthermore,
as an example for congruence in \Dfn~\ref{dfn:congruence},
the results of this algorithm can also be derived by the following left-preconditioning conversion:
\begin{eqnarray}
&&
  \tA        \PCONV \Pinv A,  \SP\SP
  \tvc{x}_k  \PCONV \vc{x}_k, \SP\SP
  \tvc{b}    \PCONV \Pinv \vc{b},   \SP\SP
  \tvc{r}_k     \PCONV \Pinv \vc{r}_k,       \SP\SP
\label{eqn:conv_precb1}
\\
&&
  \tvc{r}^\SH_0 \PCONV \PO^\T  \vc{r}^\flat_0, \SP\SP
% &&
  \tvc{p}_k \PCONV \Pinv \vc{p}_k, \SP\SP
  \tvc{u}_k \PCONV \Pinv \vc{u}_k, \SP\SP
  \tvc{q}_k \PCONV \Pinv \vc{q}_k . %, \SP\SP
 \NON
\end{eqnarray}

\subsubsection{Left-preconditioned CGS}
\label{sssec:Llprec_pcgs}

The following conversion can be used to derive another PCGS algorithm:
\begin{eqnarray}
&&
  \tA           \PCONV \Pinv A,  \SP\SP
  \tvc{x}_k     \PCONV \vc{x}_k, \SP\SP
  \tvc{b}       \PCONV \Pinv \vc{b}, \SP\SP
  \tvc{r}_k     \PCONV \vc{r}^+_k,   \SP\SP
%  \NON
\label{eqn:conv_Llprec}
\\
&&
  \tvc{r}^\SH_0 \PCONV \vc{r}^\SH_0, \SP\SP % \NON
  \tvc{p}_k \PCONV \vc{p}^+_k, \SP\SP
  \tvc{u}_k \PCONV \vc{u}^+_k, \SP\SP
  \tvc{q}_k \PCONV \vc{q}^+_k . %,\SP\SP
 \NON
\end{eqnarray}

This is the same as applying $\PO_L=\PO$ and $\PO_R=I $ to
$\tA$, $\tvc{x}_k$, and $\tvc{b}$,
that is, the left-preconditioned system.
\\

%%%%%%%%%%%%%%%%%%%%%%%%%%%%%%%%%%%%%%%%%%%%%%%%%%%%%%%%%%%%%
\refstepcounter{alg}
\label{alg:cgs_Llprec}
\noindent
\hspace{\ALGTITLEWIDTH}
{\bf Algorithm \thealg.
Left-preconditioned CGS algorithm (Left-PCGS):}
\begin{indention}{\ALGWIDTHB}
%%%%%%%%%%%%%%%%%%%%%%%%%%%%%%%%%%%%%%%%%%%%%%%%%%%%%%%%%%%%%
\noindent
$\vc{x}_0$ is the initial guess,
$\SP\vc{r}^+_0= \Pinv\L(\vc{b}-A\vc{x}_0\R),\SP\SP$
set $\SP\beta_{-1}=0,$
\\
$\LA \TL{\vc{r}}^\SH_0, \TL{\vc{r}}_0\RA
 \PCONV
 \LA \vc{r}^\SH_0, \vc{r}^+_0\RA
 \neq 0$,
e.g.,
$\vc{r}^\SH_0 = \vc{r}^+_0,\SP$
\\
For $k = 0, 1, 2, \cdots ,$ until convergence, Do:
\begin{eqnarray}
 && \vc{u}^+_k = \vc{r}^+_k + \beta_{k-1}\vc{q}^+_{k-1}, \NON
 \\
 && \vc{p}^+_k = \vc{u}^+_k + \beta_{k-1}\L(\vc{q}^+_{k-1} + \beta_{k-1}\vc{p}^+_{k-1}\R), \NON \\
 && \alpha_k = \frac
             {\LA\vc{r}^\SH_0, \vc{r}^+_k\RA}
             {\LA\vc{r}^\SH_0, \Pinv A\vc{p}^+_k\RA}, \NON
\label{eqn:Ll_cgs_alpha}
\\
 && \vc{q}^+_k = \vc{u}^+_k - \alpha_k\Pinv A\vc{p}^+_k, \NON \\
 && \vc{x}_{k+1} = \vc{x}_k + \alpha_k\L(\vc{u}^+_k + \vc{q}^+_k\R), \NON
\label{eqn:sol_Ll_pcgs}
\\
 && \vc{r}^+_{k+1} = \vc{r}^+_k - \alpha_k\Pinv A\L(\vc{u}^+_k + \vc{q}^+_k\R), \NON \\
 && \beta_k = \frac
             {\LA\vc{r}^\SH_0, \vc{r}^+_{k+1}\RA}
             {\LA\vc{r}^\SH_0, \vc{r}^+_k\RA}, \NON
\label{eqn:Ll_cgs_beta}
\end{eqnarray}
End Do
\\

\end{indention}
%%%%%%%%%%%%%%%%%%%%%%%%%%%%%%%%%%%%%%%%%%%%%%%%%%%%%%%%%%%%%

In this paper,
$\vc{r}^+_k$ denotes the residual vector
under the left-preconditioned system%
\footnote{%
The notation $\vc{r}^+_k$ will be discussed
in Sections
\ref{sssec:comparison_conventional_pcgs},
\ref{ssec:analysis_4pcgs} and \ref{sec:congruent_and_direction_pcgs}.
}%
,
its internal structure is
$\vc{r}^+_k \equiv \Pinv\vc{r}_k$,
and
this is the definition of $\vc{r}^+_k$.
Note that $\vc{p}^+_k$, $\vc{u}^+_k$, and $\vc{q}^+_k$ achieve the same purpose.
Here, $\vc{r}^+_k$ in Algorithm \ref{alg:cgs_Llprec}
provides different information to the residual vector $\vc{r}_k = \vc{b}-A\vc{x}_k$,
and the stopping criterion is
\begin{eqnarray}
\frac{\|\vc{r}^+_{k+1}\|}{\|\Pinv\vc{b}\|} \leq \varepsilon .
\label{eqn:conv_judge2}
\end{eqnarray}
Note that this is also different from (\ref{eqn:conv_judge}), 
and this is an example of incomplete judging,
because $\vc{r}^+_{k+1}$ never provides important information
about $\vc{b}-A\vc{x}_k$.
It may be thought that this is a minor issue,
but in a previous paper, we observed
that the left-preconditioned system can result in a serious problem
(see \cite{itoh2010e}, and \App~\ref{appsec:sesna_pcgs}).

This algorithm can also be derived by the following conversion:
\begin{eqnarray}
&&
  \tA           \PCONV \Pinv_L A\Pinv_R,\SP\SP
  \tvc{x}_k     \PCONV \P_R \vc{x}_k,\SP\SP
  \tvc{b}       \PCONV \Pinv_L \vc{b}, \SP\SP
  \tvc{r}_k     \PCONV \PO_R\vc{r}^+_k, \SP\SP
%  \NON
\label{eqn:sonneveld_precb_complete}
\\
&&  
  \tvc{r}^\SH_0 \PCONV \Pinvt_R\vc{r}^\SH_0, \SP\SP % \NON
  \tvc{p}_k     \PCONV \PO_R\vc{p}^+_k,\SP\SP
  \tvc{u}_k     \PCONV \PO_R\vc{u}^+_k,\SP\SP
  \tvc{q}_k     \PCONV \PO_R\vc{q}^+_k . %,\SP\SP
 \NON
\end{eqnarray}

If
 $\PO_L=\PO$ and $\PO_R=I $ are substituted into
(\ref{eqn:sonneveld_precb_complete}),
then (\ref{eqn:conv_Llprec}) is obtained.

\subsubsection{Comparison between two typical PCGS algorithms}
\label{sssec:comparison_conventional_pcgs}

Here, we compare the conventional PCGS
with the left-PCGS. We will
focus on the mathematical structures 
with respect to and in connection with their Krylov subspaces and their solution vectors.
Now,
we define a polynomial ${\it\Phi}_k(\pZ)$ symbolically
as
${\it\Phi}_k(\pZ) \equiv \alpha^{\rm PCGS}_k(\pP_k(\pZ)\pR_k(\pZ) + \pP_k(\pZ)\pR_{k+1}(\pZ))$
using Eqs.~(\ref{eqn:prec_res_pol}) through (\ref{eqn:pol_pcgs}),
then
\begin{eqnarray}
 {\it\Phi}_k(\tA)\tvc{r}_0 \equiv \alpha^{\rm PCGS}_k ( \tvc{u}_k + \tvc{q}_k )
\label{eqn:pol_phi}
\end{eqnarray}
of (\ref{eqn:x_pcgs}) in Algorithm \ref{alg:pcgs_simple},
where $k$ on both sides indicates the $k$-th iteration.
Then Eq. (\ref{eqn:ksp_cgs}) on the relation
among the structure of the solution vector, the polynomial and the Krylov subspace
can be expressed as
\begin{eqnarray}
\tvc{x}_k = \tvc{x}_\km + {\it\Phi}_\km(\tA)\tvc{r}_0
\; \in \;
\tvc{x}_0 + {\cal K}_{2k}\L(\tA, \tvc{r}_0\R).
\label{eqn:pol_cgs}
\end{eqnarray}

The conventional PCGS (Algorithm \ref{alg:cgs_precb}) is
the right-preconditioned system,
i.e.,
$\L(A\Pinv\R)\L(\PO\vc{x}\R)=\vc{b}$, and
$\vc{r}_k=\vc{b}-\L(A\Pinv\R)\L(\PO\vc{x}_k\R)$,
because this algorithm can be derived by
the right-preconditioning conversion (\ref{eqn:conv_precb2})
and this satisfies the right-preconditioned system in \Dfn~\ref{dfn:dir_prec}.
Then the vectors in \Alg~\ref{alg:cgs_precb}
and each polynomial can be detailed using (\ref{eqn:pol_pcgs}) as follows:
\begin{eqnarray}
&& \vc{r}^{\rm R}_k = \{R^{\rm R}_k(A\Pinv)\}^2\vc{r}_0, \;\;\;\; \;\;\;\;
   \vc{p}^{\rm R}_k = \{P^{\rm R}_k(A\Pinv)\}^2\vc{r}_0, %\;
\label{eqn:pol_pcgsR}
\\
&& \vc{u}^{\rm R}_k =  P^{\rm R}_k(A\Pinv) R^{\rm R}_k(A\Pinv)\vc{r}_0, \;\;\;\;
   \vc{q}^{\rm R}_k =  P^{\rm R}_k(A\Pinv) R^{\rm R}_{k+1}(A\Pinv)\vc{r}_0.
\NON
\end{eqnarray}
Here,
the vectors and polynomials of the right-preconditioned system are denoted
by superscript {\rm R},
this is also true for the scalar parameters $\alpha_k$ and $\beta_k$.
The relation between
 the solution vector
and
 its polynomial
${\it\Phi}_\km(\tA)\tvc{r}_0$
is
\begin{eqnarray}
&&\PO\vc{x}_k = \PO\vc{x}_\km + {\it\Phi}^{\rm R}_\km(A\Pinv)\vc{r}_0 , %\NON
\label{eqn:Rprec_ksp}
\end{eqnarray}
where
${\it\Phi}^{\rm R}_\km(A\Pinv)\vc{r}_0
 = \alpha^{\rm R}_\km(\vc{u}^{\rm R}_\km + \vc{q}^{\rm R}_\km)$.
This means that
 the polynomial
${\it\Phi}^{\rm R}_\km(A\Pinv)\vc{r}_0$
 organizes the solution vector as $\PO\vc{x}_k$,
 not $\vc{x}_k$ directly,
 but $\vc{x}_k$ is calculated with corrections, 
 as shown in (\ref{eqn:sol_conv_pcgs}) of \Alg~\ref{alg:cgs_precb}.

The left-PCGS (\Alg~\ref{alg:cgs_Llprec}) is the left-preconditioned system,
i.e.,
$\Pinv A\vc{x}=\Pinv\vc{b}$, and
$\vc{r}^+_k=\Pinv\L(\vc{b}-A\vc{x}_k\R)$,
because this algorithm can be derived by
the left-preconditioning conversion (\ref{eqn:conv_Llprec})
and satisfies the left-preconditioned system in \Dfn~\ref{dfn:dir_prec}.
Then the vectors in \Alg~\ref{alg:cgs_Llprec}
and each polynomial can be detailed using (\ref{eqn:pol_pcgs}) as follows:
\begin{eqnarray}
&& \vc{r}^{\rm L}_k = \{R^{\rm L}_k(\Pinv A)\}^2\vc{r}^+_0, \;\;\;\; \;\;\;\;
   \vc{p}^{\rm L}_k = \{P^{\rm L}_k(\Pinv A)\}^2\vc{r}^+_0, %\;
\label{eqn:pol_pcgsR}
\\
&& \vc{u}^{\rm L}_k =  P^{\rm L}_k(\Pinv A) R^{\rm L}_k(\Pinv A)\vc{r}^+_0, \;\;\;\;
   \vc{q}^{\rm L}_k =  P^{\rm L}_k(\Pinv A) R^{\rm L}_{k+1}(\Pinv A)\vc{r}^+_0.
\NON
\end{eqnarray}
Here, the vectors and polynomials of the left-preconditioned system are denoted
by superscript {\rm L},
this is also true for the scalar parameters $\alpha_k$ and $\beta_k$.
The relation between
 the solution vector
and
 its polynomial
${\it\Phi}_\km(\tA)\tvc{r}_0$
is
\begin{eqnarray}
&&\vc{x}_k = \vc{x}_\km + {\it\Phi}^{\rm L}_\km(\Pinv A)\vc{r}^+_0 ,
\label{eqn:sol_for_ksp_of_Llprec}
\end{eqnarray}
where
${\it\Phi}^{\rm L}_\km(\Pinv A)\vc{r}^+_0
 = \alpha^{\rm L}_\km(\vc{u}^{\rm L}_\km + \vc{q}^{\rm L}_\km)$.
Therefore,
 the polynomial
${\it\Phi}^{\rm L}_\km(\Pinv A)\vc{r}^+_0$
 organizes the solution vector directly as $\vc{x}_k$ (\Alg~\ref{alg:cgs_Llprec}).

These are summarized in \Tab~\ref{tab:summary_conv_Llprec_pcgs}.

\def\ST{\shortstack}
\renewcommand{\arraystretch}{1.2}
\begin{table}[ht]
\caption{Summary of two typical PCGS algorithms.}
% \vspace{-1mm}
\vspace{-3mm}
% \vspace{-20pt}
\begin{center}
\begin{tabular}{l|c|c}
\hline
\hline
&
\begin{tabular}{l}
Structure of\\
residual vector
\end{tabular}
&
\begin{tabular}{l}
Structure of 
the solution vector\\
for each polynomial
\end{tabular}
\\
\hline
{\makebox[33mm]{\ST{{\small Conventional (Alg.~\ref{alg:cgs_precb}) }}}}
&$\vc{r}_k = \vc{b} - (A\Pinv)(\PO\vc{x}_k)$
&
$\PO\vc{x}_k = \PO\vc{x}_\km + {\it\Phi}^{\rm R}_\km(A\Pinv)\vc{r}_0$
\\
\hline
{\makebox[33mm]{\ST{{\small Left-PCGS (Alg.~\ref{alg:cgs_Llprec})}}}}
&
$\vc{r}^+_k=\Pinv\L(\vc{b}-A\vc{x}_k\R)$
&
$\vc{x}_k = \vc{x}_\km + {\it\Phi}^{\rm L}_\km(\Pinv A)\vc{r}^+_0$
\\
\hline
\end{tabular}
\end{center}
\label{tab:summary_conv_Llprec_pcgs}
\end{table}

\subsection{Improved preconditioned CGS algorithms}
\label{ssec:improved_pcgs}

An improved PCGS algorithm has been proposed \cite{itoh2015a}.
This algorithm retains some mathematical properties that are associated
with the CGS derivation from the BiCG method
under a non-preconditioned system.
The improved PCGS algorithm from \cite{itoh2015a}
will be referred to as ``Improved1.''
Another improved PCGS algorithm will be presented,
and it will be referred to as ``Improved2.''
We note that Improved2 is mathematically equivalent to Improved1.
The stopping criterion for both algorithms is (\ref{eqn:conv_judge}).

\subsubsection{Improved1 PCGS algorithm (Improved1) \cite{itoh2015a}}
\label{sssec:improved1_pcgs}

Improved1 can be derived from the following conversion:
\begin{eqnarray}
&&
  \tA \PCONV \Pinv_LA\Pinv_R,\SP\SP
  \tvc{x}_k \PCONV \PO_R\vc{x}_k,\SP\SP
  \tvc{b}   \PCONV \Pinv_L\vc{b}, \SP\SP
  \tvc{r}_k \PCONV \Pinv_L\vc{r}_k, \SP
\label{eqn:impr1_precb_conv}
\\
&&
  \tvc{r}^\SH_0 \PCONV \Pinvt_R\vc{r}^\SH_0, \SP\SP
  \tvc{p}_k \PCONV \PO_R\vc{p}^+_k,\SP\SP
  \tvc{u}_k \PCONV \PO_R\vc{u}^+_k,\SP\SP
  \tvc{q}_k \PCONV \PO_R\vc{q}^+_k . %,\SP\SP
 \NON %\\
\end{eqnarray}

%%%%%%%%%%%%%%%%%%%%%%%%%%%%%%%%%%%%%%%%%%%%%%%%%%%%%%%%%%%%%
\refstepcounter{alg}
\label{alg:improved1_cgs_precb}
\noindent
\hspace{\ALGTITLEWIDTH}
{\bf Algorithm \thealg.
Improved PCGS algorithm (Improved1):}
\begin{indention}{\ALGWIDTHB}
%%%%%%%%%%%%%%%%%%%%%%%%%%%%%%%%%%%%%%%%%%%%%%%%%%%%%%%%%%%%%
\noindent
$\vc{x}_0$ is the initial guess,
$\SP\vc{r}_0= \vc{b}-A\vc{x}_0,\SP\SP$
set $\SP\beta_{-1}=0,$
\\
$\LA \TL{\vc{r}}^\SH_0, \TL{\vc{r}}_0\RA
 \PCONV \LA \vc{r}^\SH_0, \Pinv\vc{r}_0\RA
 \neq 0,\SP\SP$
e.g.,
$\vc{r}^\SH_0 = \Pinv\vc{r}_0,\SP$
 \\
For $k = 0, 1, 2, \cdots ,$ until convergence, Do:
\begin{eqnarray}
 && \vc{u}^+_k = \Pinv\vc{r}_k + \beta_{k-1}\vc{q}^+_{k-1}, \NON \\
 && \vc{p}^+_k = \vc{u}^+_k + \beta_{k-1}\L(\vc{q}^+_{k-1} + \beta_{k-1}\vc{p}^+_{k-1}\R), \NON \\
 && \alpha_k = \frac
 {\LA\vc{r}^\SH_0, \Pinv\vc{r}_k\RA}
 {\LA\vc{r}^\SH_0, \Pinv A\vc{p}^+_k\RA}, \NON \\
 && \vc{q}^+_k = \vc{u}^+_k - \alpha_k\Pinv A\vc{p}^+_k, \NON \\
 && \vc{x}_{k+1} = \vc{x}_k + \alpha_k\L(\vc{u}^+_k + \vc{q}^+_k\R), \NON \\
 && \vc{r}_{k+1} = \vc{r}_k - \alpha_kA\L(\vc{u}^+_k + \vc{q}^+_k\R), \NON \\
 && \beta_k = \frac
 {\LA\vc{r}^\SH_0, \Pinv\vc{r}_{k+1}\RA}
 {\LA\vc{r}^\SH_0, \Pinv\vc{r}_k\RA}, \NON
\end{eqnarray}
End Do
% \\

\end{indention}
%%%%%%%%%%%%%%%%%%%%%%%%%%%%%%%%%%%%%%%%%%%%%%%%%%%%%%%%%%%%%

\subsubsection{Improved2 PCGS algorithm (Improved2)}
\label{sssec:improved2_pcgs}

Improved2 can be derived from the following conversion:
\begin{eqnarray}
&&
  \tA       \PCONV \Pinv_L A \Pinv_R,  \SP\SP
  \tvc{x}_k \PCONV \PO_R      \vc{x}_k, \SP\SP
  \tvc{b}   \PCONV \Pinv_L   \vc{b}, \SP\SP
  \tvc{r}_k     \PCONV \Pinv_L  \vc{r}_k, \SP
\label{eqn:impr2_precb_conv}
\\
&&
  \tvc{r}^\SH_0 \PCONV \Pinvt_R \vc{r}^\SH_0, \SP\SP
  \tvc{p}_k \PCONV \Pinv_L \vc{p}_k, \SP\SP
  \tvc{u}_k \PCONV \Pinv_L \vc{u}_k, \SP\SP
  \tvc{q}_k \PCONV \Pinv_L \vc{q}_k . %,\SP\SP 
 \NON %\\
\end{eqnarray}

Note that this conversion is different than (\ref{eqn:impr1_precb_conv}) for $\tvc{p}_k$, $\tvc{u}_k$, and $\tvc{q}_k$.
\\

%%%%%%%%%%%%%%%%%%%%%%%%%%%%%%%%%%%%%%%%%%%%%%%%%%%%%%%%%%%%%
\refstepcounter{alg}
\label{alg:improved2_cgs_precb}
\noindent
\hspace{\ALGTITLEWIDTH}
{\bf Algorithm \thealg.
Another improved PCGS algorithm (Improved2):}
\begin{indention}{\ALGWIDTHB}
%%%%%%%%%%%%%%%%%%%%%%%%%%%%%%%%%%%%%%%%%%%%%%%%%%%%%%%%%%%%%
\noindent
$\vc{x}_0$ is the initial guess,
$\SP\vc{r}_0= \vc{b}-A\vc{x}_0,\SP\SP$
set $\SP\beta_{-1}=0,$
\\
$\LA \TL{\vc{r}}^\SH_0, \TL{\vc{r}}_0\RA
 \PCONV \LA \vc{r}^\SH_0, \Pinv\vc{r}_0\RA
 \neq 0,\SP\SP$
e.g.,
$\vc{r}^\SH_0 = \Pinv\vc{r}_0$,
\\
For $k = 0, 1, 2, \cdots ,$ until convergence, Do:
\begin{eqnarray}
 && \vc{u}_k = \vc{r}_k + \beta_{k-1}\vc{q}_{k-1}, \NON \\
 && \vc{p}_k = \vc{u}_k + \beta_{k-1}\L(\vc{q}_{k-1} + \beta_{k-1}\vc{p}_{k-1}\R), \NON
\\
 && \alpha_k = \frac
              {\LA \Pinvt\vc{r}^\SH_0, \vc{r}_k\RA}
              {\LA \Pinvt\vc{r}^\SH_0, A\Pinv\vc{p}_k\RA}, \NON
\label{eqn:impr2_cgs_alpha}
\\
 && \vc{q}_k = \vc{u}_k - \alpha_kA\Pinv\vc{p}_k, \NON
\\
 && \vc{x}_{k+1} = \vc{x}_k + \alpha_k\Pinv\L(\vc{u}_k + \vc{q}_k\R), \NON
\\
 && \vc{r}_{k+1} = \vc{r}_k - \alpha_kA\Pinv\L(\vc{u}_k + \vc{q}_k\R), \NON \\
 && \beta_k = \frac
              {\LA \Pinvt\vc{r}^\SH_0, \vc{r}_{k+1}\RA}
              {\LA \Pinvt\vc{r}^\SH_0, \vc{r}_k\RA}, \NON
\label{eqn:impr2_cgs_beta}
\end{eqnarray}
End Do
\\

\end{indention}
%%%%%%%%%%%%%%%%%%%%%%%%%%%%%%%%%%%%%%%%%%%%%%%%%%%%%%%%%%%%%

\subsection{Analysis of the four kinds of PCGS algorithms}
\label{ssec:analysis_4pcgs}

We will now analyze and compare the four PCGS algorithms presented above.

\begin{figure}[htbp]
\begin{center}
\resizebox*{\ALGCHARTSIZE\columnwidth}{!}{
\includegraphics*{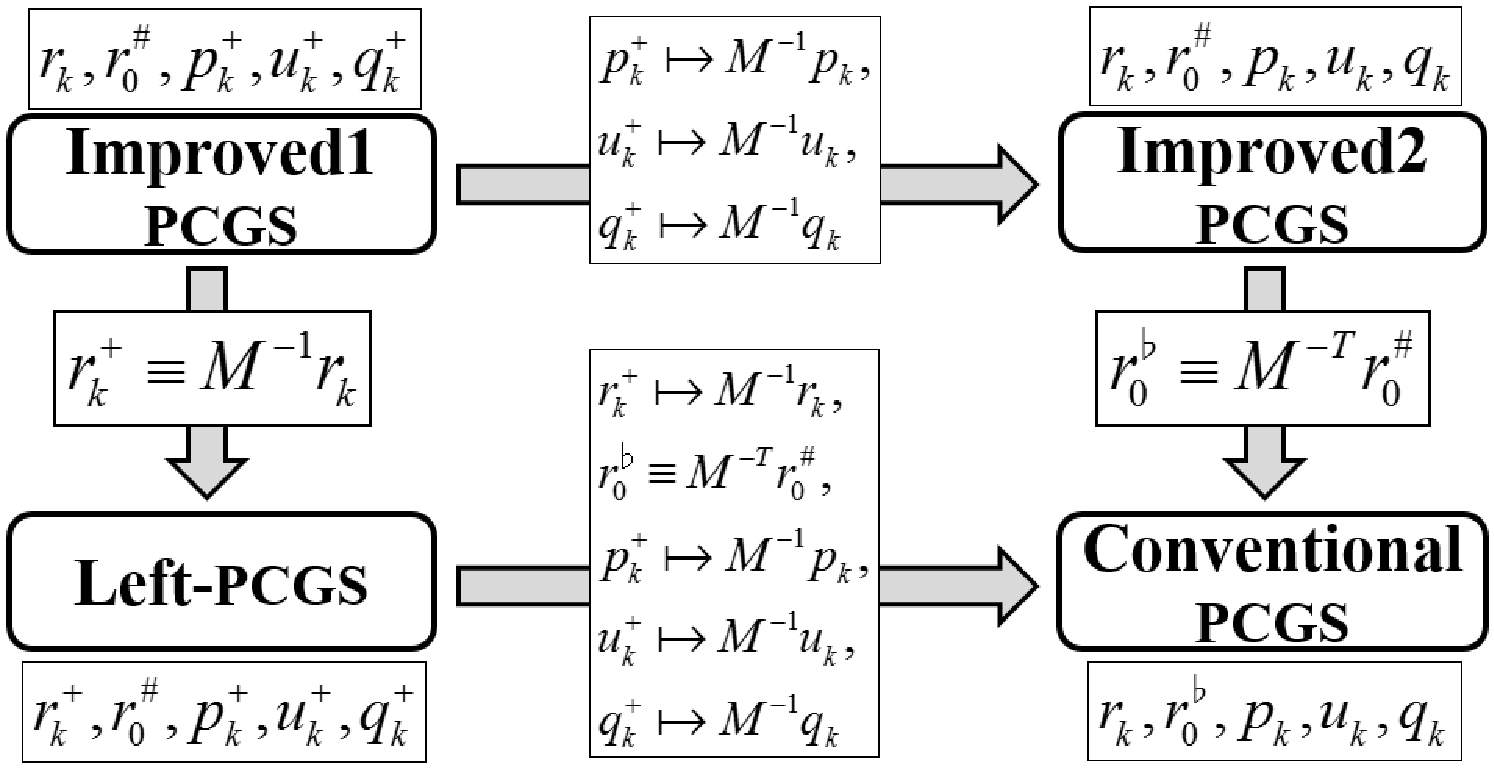}}
\caption{Relations between the four different PCGS algorithms.
$\mapsto$ : Splitting left vector to right members (preconditioner and vector),
$\equiv$ : Substituting left vector for right members.}
\label{fig:PCGS_ConvLlprecImproved}
\end{center}
\vspace{-13pt}
\end{figure}

We split the residual vector of the left-PCGS (\Alg~\ref{alg:cgs_Llprec})
$\vc{r}^+_k$ into
\begin{eqnarray}
&&\vc{r}^+_k \mapsto \Pinv\vc{r}_k, \SP\SP (k = 0, 1, 2, \cdots)
\label{eqn:r_plus_convert}
\end{eqnarray}
and give the necessary deformations;
then,
the left-PCGS (\Alg~\ref{alg:cgs_Llprec}) is reduced to
Improved1 (\Alg~\ref{alg:improved1_cgs_precb}).
Alternatively, we can derive
\Alg~\ref{alg:cgs_Llprec} from \Alg~\ref{alg:improved1_cgs_precb}
by substituting $\Pinv\vc{r}_k$ for $\vc{r}^+_k$,
that is, $\vc{r}^+_k \equiv \Pinv\vc{r}_k$.
By this means,
we can explain the relationships between the four kinds of PCGS algorithms,
as shown in \Fig~\ref{fig:PCGS_ConvLlprecImproved}.

In addition,
if we apply (\ref{eqn:r_plus_convert}) to 
(\ref{eqn:sol_for_ksp_of_Llprec})
to obtain the structure of the polynomial ${\it\Phi}_k(\tA)\tvc{r}_0$ of \Alg~\ref{alg:cgs_Llprec},
then
\begin{eqnarray}
&&
{\it\Phi}^{\rm L}_k(\Pinv A)\vc{r}^+_0
\SP\SP \mapsto \SP\SP
{\it\Phi}^{\rm L}_k(\Pinv A)\Pinv\vc{r}_0
=
\Pinv{\it\Phi}^{\rm L}_k(A\Pinv)\vc{r}_0. \NON
\end{eqnarray}
Therefore,
the system of Improved1 (\Alg~\ref{alg:improved1_cgs_precb})
is coordinative to that of the left-PCGS (\Alg~\ref{alg:cgs_Llprec}).
Furthermore,
the structure of the solution vector and the polynomial is
\begin{eqnarray}
&&
\vc{x}_k = \vc{x}_\km + {\it\Phi}^{\rm L}_\km(\Pinv A)\vc{r}^+_0
\SP\SP \mapsto \SP\SP
\vc{x}_k = \vc{x}_\km + \Pinv{\it\Phi}^{\rm L}_\km(A\Pinv)\vc{r}_0 .
\label{eqn:structure_of_solution}
\end{eqnarray}
Here,
the polynomial is denoted as ${\it\Phi}^{\rm L}_\km$, 
despite the preconditioned matrix being $A\Pinv$ in the right-hand side
of (\ref{eqn:structure_of_solution}),
because
the direction of a preconditioned system
is different from
the direction of a preconditioning conversion
(see Section \ref{sec:congruent_and_direction_pcgs}).

Improved2 (\Alg~\ref{alg:improved2_cgs_precb}) is equivalent to
Improved1.
Both algorithms have important advantages over the left-PCGS,
because their residual vector is $\vc{r}_k$,
and their stopping criterion is (\ref{eqn:conv_judge}),
not $\|\vc{r}^+_{k+1}\|/\|\Pinv\vc{b}\|$.

\Tab\ref{tab:summary_All_pcgs} shows
the structure of the residual vector
and 
the structure of the solution vector for each polynomial
of the four PCGS algorithms.

\def\ST{\shortstack}
\renewcommand{\arraystretch}{1.2}
\begin{table}[h]
\caption{Summary of the four PCGS algorithms.}
\vspace{-3mm}
% \vspace{-2mm}
\begin{center}
\begin{tabular}{l|c|c}
\hline
\hline
&
\begin{tabular}{l}
Structure of\\
residual vector
\end{tabular}
&
\begin{tabular}{l}
Structure of the
solution vector\\
for each polynomial
\end{tabular}
\\
\hline
{\makebox[33mm]{\ST{
{\small
Conventional (Alg.~\ref{alg:cgs_precb})
}
}}}
&
{\small
$\vc{r}_k = \vc{b} - (A\Pinv)(\PO\vc{x}_k)$
}
&
{\small
$\PO\vc{x}_k = \PO\vc{x}_\km + {\it\Phi}^{\rm R}_\km(A\Pinv)\vc{r}_0$
}
\\
\hline
{\makebox[33mm]{\ST{
{\small
Left-PCGS  (Alg.~\ref{alg:cgs_Llprec})
}
}}}
&
{\small
$\vc{r}^+_k=\Pinv\L(\vc{b}-A\vc{x}_k\R)$
}
&
{\small
$\vc{x}_k = \vc{x}_\km + {\it\Phi}^{\rm L}_\km(\Pinv A)\vc{r}^+_0$
}
\\
\hline
{\makebox[33mm]{\ST{
{\small
Improved1 
(Alg.~\ref{alg:improved1_cgs_precb})
}
}}}
&
\multicolumn{1}{l|}
{\small
$\Pinv\vc{r}_k = $
}
&
\\
\cline{1-1}
{\makebox[33mm]{\ST{
{\small
Improved2 
(Alg.~\ref{alg:improved2_cgs_precb})
}
}}}
&
\multicolumn{1}{r|}
{\small
$\Pinv\L(\vc{b} - (A\Pinv)(\PO\vc{x}_k)\R)$
}
&
\multicolumn{1}{c}{\raisebox{1.5ex}[0pt]{
{\makebox[40mm]{\ST{
{\small
$\vc{x}_k = \vc{x}_\km + \Pinv{\it\Phi}^{\rm L}_\km(A\Pinv)\vc{r}_0$
}
}}}
}}
\\
\hline
\end{tabular}
\vspace{-13pt}
\end{center}
\label{tab:summary_All_pcgs}
\end{table}

In this summary, we see that
the structures of the polynomial organizing the solution vector differ:
\begin{eqnarray}
&& {\it\Phi}^{\rm R}_\km(A\Pinv)\vc{r}_0
 \neq
   {\it\Phi}^{\rm L}_\km(A\Pinv)\vc{r}_0.
\label{eqn:differ_bothKSP}
\end{eqnarray}
In other words,
  ${\it\Phi}^{\rm R}_\km(A\Pinv)\vc{r}_0$
  of the conventional PCGS (right-preconditioned) system
is different from
  ${\it\Phi}^{\rm L}_\km(A\Pinv)\vc{r}_0$
of both improved PCGS (coordinative to the left-preconditioned) systems,
because
the scalar parameters $\alpha_k$ and $\beta_k$
are not equivalent
(see Section \ref{ssec:dir_and_same_pcgs})
\cite{itoh2015a,itoh2016b}.
This also affects the polynomial ${\it\Phi}_k(\tA)\tvc{r}_0$
for the following reason.
Although,
superficially, 
the solution vector for both
  the conventional PCGS (\Alg~\ref{alg:cgs_precb})
and
 Improved2 (\Alg~\ref{alg:improved2_cgs_precb})
have the same recurrence formula,
i.e.,
$\vc{x}_{k+1} = \vc{x}_k + \alpha_k\Pinv\L(\vc{u}_k + \vc{q}_k\R)$,
each recurrence formula belongs to a different system
because
the components of the conventional PCGS are
$\alpha^{\rm R}_k$,
$\vc{u}^{\rm R}_k$, and
$\vc{q}^{\rm R}_k$,
and
those of Improved2 are
$\alpha^{\rm L}_k$,
$\vc{u}^{\rm L}_k$, and 
$\vc{q}^{\rm L}_k$.

The structure of the residual vector
of Improved1 (\Alg~\ref{alg:improved1_cgs_precb})
and Improved2 (\Alg~\ref{alg:improved2_cgs_precb})
appears as
$\Pinv\vc{r}_k=\Pinv\L(\vc{b} - (A\Pinv)(\PO\vc{x}_k)\R)$
in \Tab\ref{tab:summary_All_pcgs}
because they are from the left-PCGS,
and the structure of their polynomial is
$\Pinv{\it\Phi}^{\rm L}_\km(A\Pinv)\vc{r}_0$.

Note that
(\ref{eqn:differ_bothKSP}) will 
be confirmed in Section \ref{sec:congruent_and_direction_pcgs}
and
will be shown numerically
in Section \ref{sec:numerical_experiments}.

%%%%%%  section  %%%%%%%%%%%%%%%%%%%%%%%%%%%%%%%%%%%%%%%%%%%%%%%%%%%%%%%%%%%%%%%%%
\section{Congruence of preconditioning conversion,
and direction of preconditioned system for the CGS}
\label{sec:congruent_and_direction_pcgs}

In the previous section, we defined
the general direction of a preconditioned system for CGS (see
\Dfn~\ref{dfn:dir_prec}).
However,
the direction of a preconditioned system is different from
the direction of a preconditioning conversion.
We will show that the direction of a preconditioned system is switched
by the construction of the initial shadow residual vector (ISRV).

\subsection{Congruence of preconditioning conversion for the PCGS}
\label{ssec:congruent_pcgs}

Here,
we consider the congruence of \Dfn~\ref{dfn:congruence}
for the PCGS in the following \prop.

\begin{proposition}[Congruency]
\label{prop:congruency}
There is congruence to a PCGS algorithm
in the direction of the preconditioning conversion.
\end{proposition}

{\bf Proof }
We have already shown instances of this.
For example,
\Alg~\ref{alg:cgs_precb} can be derived
by the two-sided conversion (\ref{eqn:conv_precb}),
and
if $\PO_L=I$, $\PO_R=\PO $, and
the conversion (\ref{eqn:conv_precb}) is reduced
to (\ref{eqn:conv_precb2}),
then \Alg~\ref{alg:cgs_precb} is derived.
If $\PO_L=\PO$ and $\PO_R=I $,
that is (\ref{eqn:conv_precb1}),
then
\Alg~\ref{alg:cgs_precb} can be derived.
The other preconditioned algorithms
(Algorithms~\ref{alg:cgs_Llprec}, \ref{alg:improved1_cgs_precb},
and \ref{alg:improved2_cgs_precb})
and their corresponding preconditioning conversions are also the same.
\qquad
$\hfill\Box$
\\

Although this property has been repeatedly discussed in the literature,
it should be considered when evaluating
the direction of a preconditioned system.

\subsection{Direction of a preconditioned system and that of the PCGS}
\label{ssec:dir_and_same_pcgs}

The direction of a preconditioned system is different from
the direction of a preconditioning conversion.

\begin{proposition}
\label{prop:pcgs_biortho_biconj}
The direction of a preconditioned system %for a CGS algorithm
 is determined by
the operations of $\alpha_k$ and $\beta_k$ in each PCGS algorithm.
These intrinsic operations are based on
 biorthogonality $(\tvc{r}^\SH_0, \tvc{r}_k)$
and
 biconjugacy $(\tvc{r}^\SH_0, \tA\tvc{p}_k)$.
\end{proposition}

{\bf Proof }
The operations of biorthogonality and biconjugacy
in each PCGS algorithm
and
the structure of the solution vector for each polynomial
are shown below.
The underlined inner products are the actual descriptions for each PCGS algorithm.

Only the conventional PCGS (\Alg~\ref{alg:cgs_precb}) algorithm has
the ISRV in the form $\vc{r}^\flat_0$; in all other algorithms, it is $\vc{r}^\SH_0$.
The ISRV $\vc{r}^\flat_0$ never splits into $\Pinvt\vc{r}^\SH_0$ in this algorithm,
and the preconditioned coefficient matrix for the biconjugacy
is fixed as $A\Pinv$,
that is, the right-preconditioned system.

\begin{itemize}

\item
Conventional (\Alg~\ref{alg:cgs_precb}) :
\begin{eqnarray}
\vc{r}^\flat_0 &=& \vc{r}_0, \NON
\\
\L(\tvc{r}^\SH_0, \tvc{r}_k\R)
 &\PCONV& \L(\PO^\T_L\vc{r}^\flat_0,\SP \Pinv_L\vc{r}_k\R)
 =        \UL{ \L(\vc{r}^\flat_0,\SP \vc{r}_k\R) }, \NON
\\
\L(\tvc{r}^\SH_0, \tA\tvc{p}_k\R)
 &\PCONV&      \L(\PO^\T_L\vc{r}^\flat_0,\SP (\Pinv_L A\Pinv_R)(\Pinv_L\vc{p}_k)\R)
 =             \UL{ \L(\vc{r}^\flat_0,\SP (A\Pinv ) \vc{p}_k\R) }, \NON
\\
\PO\vc{x}_k &=& \PO\vc{x}_\km + {\it\Phi}^{\rm R}_\km(A\Pinv)\vc{r}_0. \NON
\end{eqnarray}

\item
Left-PCGS (\Alg~\ref{alg:cgs_Llprec}) :
\begin{eqnarray}
\vc{r}^\SH_0 &=& \vc{r}^+_0, \NON
\\
\L(\tvc{r}^\SH_0, \tvc{r}_k\R)
 &\PCONV&  \UL{ \L(\vc{r}^\SH_0,\SP \vc{r}^+_k\R) }, \NON
\\
\L(\tvc{r}^\SH_0, \tA\tvc{p}_k\R)
 &\PCONV&  \UL{ \L(\vc{r}^\SH_0,\SP (\Pinv A)\vc{p}^+_k\R) }, \NON
\\
\vc{x}_k &=& \vc{x}_\km + {\it\Phi}^{\rm L}_\km(\Pinv A)\vc{r}^+_0. \NON
\mfill{0.37}
\end{eqnarray}

\item
Improved1 (\Alg~\ref{alg:improved1_cgs_precb}) :
\begin{eqnarray}
\vc{r}^\SH_0 &=& \Pinv\vc{r}_0, \NON
\\
\L(\tvc{r}^\SH_0, \tvc{r}_k\R)
 &\PCONV& \L(\Pinvt_R\vc{r}^\SH_0,\SP \Pinv_L\vc{r}_k\R)
 =        \UL{ \L(\vc{r}^\SH_0,\SP \Pinv\vc{r}_k\R) }, \NON
\\
\L(\tvc{r}^\SH_0, \tA\tvc{p}_k\R)
 &\PCONV&  \L(\Pinvt_R\vc{r}^\SH_0,\SP (\Pinv_L A\Pinv_R)(\PO_R\vc{p}^+_k)\R)
 =         \UL{ \L(\vc{r}^\SH_0, (\Pinv A) \vc{p}^+_k \R) }, \NON
\\
\vc{x}_k &=& \vc{x}_\km + \Pinv{\it\Phi}^{\rm L}_\km(A\Pinv)\vc{r}_0. \NON
\end{eqnarray}

\item
Improved2 (\Alg~\ref{alg:improved2_cgs_precb}) :
\begin{eqnarray}
\vc{r}^\SH_0 &=& \Pinv\vc{r}_0, \NON
\\
\L(\tvc{r}^\SH_0, \tvc{r}_k\R)
 &\PCONV& \L(\Pinvt_R\vc{r}^\SH_0,\SP \Pinv_L\vc{r}_k\R) \NON
  =
\UL{ \L(\Pinvt\vc{r}^\SH_0,\SP \vc{r}_k\R) }
 =         \L(\vc{r}^\SH_0,\SP \Pinv\vc{r}_k\R) , \NON
\\
\L(\tvc{r}^\SH_0, \tA\tvc{p}_k\R)
 &\PCONV&   \L(\Pinvt_R\vc{r}^\SH_0,\SP (\Pinv_L A\Pinv_R)(\Pinv_L\vc{p}_k)\R) \NON
\\
 &=&        \UL{ \L(\Pinvt\vc{r}^\SH_0,\SP A(\Pinv \vc{p}_k)\R) }
  =              \L(\vc{r}^\SH_0,\SP (\Pinv A)(\Pinv \vc{p}_k)\R), \NON
\\
\vc{x}_k &=& \vc{x}_\km + \Pinv{\it\Phi}^{\rm L}_\km(A\Pinv)\vc{r}_0. \NON
\qquad
\hfill\Box
\end{eqnarray}
\end{itemize}

We present the following proposition and \cor.

\begin{proposition}
\label{prop:pcgs_biortho}
On the structure of biorthogonality $(\tvc{r}^\SH_0, \tvc{r}_k)$
in the iterated part of each PCGS algorithm,
there exists a single preconditioning operator
between
  $\vc{r}_k$ (basic form of the residual vector)
and
  $\vc{r}^\SH_0$ (basic form of the ISRV) such that
  $\Pinv$ operates on $\vc{r}_k$
or
  $\Pinvt$ operates on $\vc{r}^\SH_0$.
\end{proposition}

{\bf Proof }
We split
$\vc{r}^\flat_0 \mapsto \Pinvt\vc{r}^\SH_0$ and
$\vc{r}^+_k     \mapsto \Pinv \vc{r}_k$
in Algorithms \ref{alg:cgs_precb} to \ref{alg:improved2_cgs_precb},
and obtain
\begin{eqnarray}
\L( \tvc{r}^\SH_0, \tvc{r}_k \R)
   &=&   \UL{ \L( \vc{r}^\flat_0, \SP \vc{r}_k \R) }
 \mapsto \L( \Pinvt \vc{r}^\SH_0, \SP \vc{r}_k \R) , \NON
\label{eqn:cgs_precb_biortho_detail}
\\
\L( \tvc{r}^\SH_0, \tvc{r}_k \R)
   &=&   \UL{ \L(\vc{r}^\SH_0, \SP \vc{r}^+_k\R) }
 \mapsto \L( \vc{r}^\SH_0, \SP \Pinv \vc{r}_k \R) , \NON
\label{eqn:cgs_Llprec_biortho_detail}
\\
\L( \tvc{r}^\SH_0, \tvc{r}_k \R)
   &=&   \UL{ \L(\vc{r}^\SH_0, \SP \Pinv\vc{r}_k\R) }, \NON
\label{eqn:cgs_impr1_biortho_detail}
\\
\L( \tvc{r}^\SH_0, \tvc{r}_k \R)
   &=&   \UL{ \L(\Pinvt\vc{r}^\SH_0, \SP \vc{r}_k\R) }
    =    \L( \vc{r}^\SH_0, \SP \Pinv\vc{r}_k \R) . \NON
\label{eqn:cgs_impr2_biortho_detail}
\end{eqnarray}
The underlined inner products are the actual descriptions for each PCGS algorithm.

In addition, for
the two-sided conversion,
we obtain
\begin{eqnarray}
\L( \tvc{r}^\SH_0, \tvc{r}_k \R)
   &=&   \L(\Pinvt_R\vc{r}^\SH_0, \SP \Pinv_L\vc{r}_k\R)
    =    \L(\Pinvt\vc{r}^\SH_0, \SP \vc{r}_k\R)
    =    \L( \vc{r}^\SH_0, \SP \Pinv\vc{r}_k \R) . \NON
\label{eqn:cgs_two_biortho_detail}
\qquad
\hfill\Box
\end{eqnarray}

\begin{corollary}
\label{cor:pcgs_biconj}
On the structure of biconjugacy $(\tvc{r}^\SH_0, \tA \tvc{p}_k)$
in the iterated part of each PCGS algorithm,
there exists a single preconditioning operator
between
  $A$ (coefficient matrix)
and
  $\vc{r}^\SH_0$ (basic form of the ISRV),
such that
  $\Pinv$ operates on $A$
or
  $\Pinvt$ operates on $\vc{r}^\SH_0$.
Furthermore,
there exists a single preconditioning operator
between
  $A$
and
  $\vc{p}_k$ (basic form of probing direction vector).
\end{corollary}

From 
Propositions~\ref{prop:pcgs_biortho_biconj} and \ref{prop:pcgs_biortho}
and \Cor~\ref{cor:pcgs_biconj},
the intrinsic operations on the biorthogonality and the biconjugacy
for the four PCGS algorithms
have the same matrix and vector structures,
even though the superficial descriptions of these algorithms are different.

\subsection{ISRV switches the direction of the preconditioned system for the CGS}
\label{ssec:isrv}

Although the
mathematical properties of
the conventional PCGS (\Alg~\ref{alg:cgs_precb} : right-preconditioned system)
and
Improved2 (\Alg~\ref{alg:improved2_cgs_precb}: coordinative to the left-preconditioned system)
are quite different,
the structures of these algorithms are very similar.
This can be seen by
replacing $\Pinvt\vc{r}^\SH_0$ with $\vc{r}^\flat_0$
in \Alg~\ref{alg:improved2_cgs_precb},
and
in the initial part, we have
\begin{eqnarray}
\L( \tvc{r}^\SH_0, \tvc{r}_0 \R)
 =  \L( \Pinvt_R\vc{r}^\SH_0, \Pinv_L\vc{r}_0 \R)
&=& \L( \Pinvt\vc{r}^\SH_0, \vc{r}_0 \R) \NON
\\
&\equiv& \L( \vc{r}^\flat_0, \vc{r}_0 \R) \neq 0,\SP\SP
{\rm e.g.,} \SP\SP \vc{r}^\flat_0=\vc{r}_0 ; \NON
\end{eqnarray}
then
\Alg~\ref{alg:improved2_cgs_precb} becomes \Alg~\ref{alg:cgs_precb}.

\begin{theorem}
\label{thm:isrv_switch}
The direction of a preconditioned system for the CGS is switched
by the construction and setting of the ISRV.
\end{theorem}

{\bf Proof }
\Prop~\ref{prop:pcgs_biortho_biconj} shows that
the direction of a preconditioned system for the CGS algorithm is
determined by the structures of the biorthogonality and the biconjugacy.
Here,
we show that
their structures 
 are
switched by the ISRV.
The underlined inner products are the actual operators for each PCGS algorithm.

\begin{itemize}

\item
ISRV1 : $\vc{r}^\SH_0 = \Pinv \vc{r}_0$ (Based on left conversion)
\begin{eqnarray}
\LA \tvc{r}^\SH_0, \tvc{r}_0 \RA
&=&
\LA \Pinvt_R \vc{r}^\SH_0, \Pinv_L \vc{r}_0 \RA
= \LA \vc{r}^\SH_0, \Pinv \vc{r}_0 \RA \neq 0, \NON
\\
&&{\rm e.g., }\SP\SP \vc{r}^\SH_0 = \Pinv \vc{r}_0 . \NON
\end{eqnarray}

\item
ISRV2 : $\vc{r}^\SH_0 = \PO^\T \vc{r}_0$ (Based on right conversion)
\begin{eqnarray}
\LA \tvc{r}^\SH_0, \tvc{r}_0 \RA
&=&
\LA \Pinvt_R \vc{r}^\SH_0, \Pinv_L \vc{r}_0 \RA
= \LA \Pinvt \vc{r}^\SH_0, \vc{r}_0 \RA \neq 0, \NON
\\
&&{\rm e.g., }\SP\SP \Pinvt \vc{r}^\SH_0 = \vc{r}_0 \SP
\rightarrow \vc{r}^\SH_0 = \PO^\T \vc{r}_0 . \NON
\end{eqnarray}
\end{itemize}

If we apply ISRV2 to \Alg~\ref{alg:improved2_cgs_precb},
then \Alg~\ref{alg:improved2_cgs_precb} is equivalent to \Alg~\ref{alg:cgs_precb}
with $\vc{r}^\flat_0 = \vc{r}_0$.
Here,
the operations in the iterated part 
of both algorithms are the same,
because $\vc{r}^\flat_0$ in \Alg~\ref{alg:cgs_precb} is equivalent to $\Pinvt\vc{r}^\SH_0$.
\begin{enumerate}
\item[1-1)] \Alg~\ref{alg:cgs_precb} with $\vc{r}^\flat_0 = \vc{r}_0$ as ISRV
 is the right-preconditioned system:
\begin{eqnarray}
\L(\tvc{r}^\SH_0, \tvc{r}_k\R)
 &=&  
    \UL{ \L(\vc{r}^\flat_0,\SP \vc{r}_k\R) }
 =         \L(\vc{r}_0,\SP \vc{r}_k\R) , \NON
\\
\L(\tvc{r}^\SH_0, \tA\tvc{p}_k\R)
 &=&
    \UL{ \L(\vc{r}^\flat_0,\SP A(\Pinv \vc{p}_k)\R) }
  =              \L(\vc{r}_0,\SP (A\Pinv) \vc{p}_k \R). \NON
\end{eqnarray}
\item[1-2)] \Alg~\ref{alg:improved2_cgs_precb} with ISRV2
 is equivalent to the right-preconditioned system:
\begin{eqnarray}
\L(\tvc{r}^\SH_0, \tvc{r}_k\R)
 &=&  
    \UL{ \L(\Pinvt\vc{r}^\SH_0,\SP \vc{r}_k\R) }
 =  \L(\Pinvt(\PO^\T\vc{r}_0),\SP \vc{r}_k\R)
 =         \L(\vc{r}_0,\SP \vc{r}_k\R) , \NON
\\
\L(\tvc{r}^\SH_0, \tA\tvc{p}_k\R)
 &=&
    \UL{ \L(\Pinvt\vc{r}^\SH_0,\SP A(\Pinv \vc{p}_k)\R) }
  = \L(\Pinvt(\PO^\T\vc{r}_0),\SP A\Pinv \vc{p}_k\R) \NON
\\
 &=&             \L(\vc{r}_0,\SP (A\Pinv) \vc{p}_k \R). \NON
\end{eqnarray}
\end{enumerate}

Alternatively, if we apply
$\vc{r}^\flat_0 = \Pinvt\Pinv \vc{r}_0$ (we will call this ISRV9
\footnote{
Although ISRV9 is not sequential with respect to ISRV1 and ISRV2,
the designation is consist with our research.
}%
) to
\Alg~\ref{alg:cgs_precb},
then \Alg~\ref{alg:cgs_precb} is equivalent to \Alg~\ref{alg:improved2_cgs_precb}
with ISRV1.
Here,
the operations in the iterated part 
of both algorithms are the same
as mentioned above.
\begin{enumerate}
\item[2-1)] \Alg~\ref{alg:improved2_cgs_precb} with ISRV1 is coordinative to the left-preconditioned system:
\begin{eqnarray}
\L(\tvc{r}^\SH_0, \tvc{r}_k\R)
 &=&  
\UL{ \L(\Pinvt\vc{r}^\SH_0,\SP \vc{r}_k\R) }
 =   \L(\Pinvt(\Pinv\vc{r}_0),\SP \vc{r}_k\R)
 =   \L(\Pinv\vc{r}_0,\SP \Pinv\vc{r}_k\R)
 , \NON
\\
\L(\tvc{r}^\SH_0, \tA\tvc{p}_k\R)
 &=&        \UL{ \L(\Pinvt\vc{r}^\SH_0,\SP (A\Pinv) \vc{p}_k\R) }
  =              \L(\Pinvt(\Pinv\vc{r}_0),\SP (A\Pinv) \vc{p}_k \R) \NON
\\
 &&
  \;\;\;\;\;\;\;\;
  \;\;\;\;\;\;\;
  \;\;\;\;\;\;\;
  \;\;\;\;\;\;\;
 = \L(\Pinv\vc{r}_0,\SP (\Pinv A)(\Pinv\vc{p}_k) \R)
. \NON
\end{eqnarray}
\item[2-2)] \Alg~\ref{alg:cgs_precb} with ISRV9 is coordinative to the left-preconditioned system:
\begin{eqnarray}
\L(\tvc{r}^\SH_0, \tvc{r}_k\R)
 &=&  
\UL{ \L(\vc{r}^\flat_0,\SP \vc{r}_k\R) }
 =   \L(\Pinvt(\Pinv\vc{r}_0),\SP \vc{r}_k\R)
 =   \L(\Pinv\vc{r}_0,\SP \Pinv\vc{r}_k\R)
 , \NON
\\
\L(\tvc{r}^\SH_0, \tA\tvc{p}_k\R)
 &=&        \UL{ \L(\vc{r}^\flat_0,\SP (A\Pinv) \vc{p}_k\R) }
  =              \L(\Pinvt(\Pinv\vc{r}_0),\SP (A\Pinv) \vc{p}_k \R) \NON
\\
 &&
  \;\;\;\;\;\;\;\;
  \;\;\;\;\;\;\;
  \;\;\;\;\;\;\;
  \;\;\;\;\;\;\;
 = \L(\Pinv\vc{r}_0,\SP (\Pinv A)(\Pinv\vc{p}_k) \R)
. \NON
\qquad
\hfill\Box
\end{eqnarray}
\end{enumerate}

If we change Improved2 (\Alg~\ref{alg:improved2_cgs_precb})
 to
Improved1 (\Alg~\ref{alg:improved1_cgs_precb}),
then we will obtain the same results.

In this section,
we have confirmed (\ref{eqn:differ_bothKSP}),
i.e.,
${\it\Phi}^{\rm R}_\km(A\Pinv)\vc{r}_0
 \neq
 {\it\Phi}^{\rm L}_\km(A\Pinv)\vc{r}_0$,
by \Prop~\ref{prop:pcgs_biortho_biconj}
with
\Thm~\ref{thm:isrv_switch}.
In the next section,
\Thm~\ref{thm:isrv_switch} is verified numerically.

%%%%%%  section  %%%%%%%%%%%%%%%%%%%%%%%%%%%%%%%%%%%%%%%%%%%%%%%%%%%%%%%%%%%%%%%%%
\section{Numerical experiments}
\label{sec:numerical_experiments}

Convergence of
the four PCGS algorithms of Section~\ref{sec:pcgs_alg}
is
confirmed in Section~\ref{ssec:num_exp1} by evaluating three cases.
Furthermore, in Section~\ref{ssec:num_exp2}, the ability of the ISRV to
switch the direction of the preconditioned system (as discussed
in Section~\ref{ssec:isrv}) is verified, as well as \Thm~\ref{thm:isrv_switch}.

\subsection{Comparison of the four PCGS algorithms}
\label{ssec:num_exp1}

The test problems were generated
by building real nonsymmetric matrices corresponding to linear systems
taken from the University of Florida Sparse Matrix Collection \cite{ufl}
and the Matrix Market \cite{matrixmarket}.
The RHS vector $\vc{b}$ of (\ref{eqn:linear})
was generated by setting all elements of the exact solution vector
$\vc{x}_{\rm exact}$
to 1.0 and substituting this into (\ref{eqn:linear}).
The solution algorithm was implemented using the sequential mode
of the Lis numerical computation library (version 1.1.2 \cite{lis})
in double precision,
with the compiler options registered in the Lis ``Makefile.''
Furthermore,
we set the initial solution to $\vc{x}_0 = \vc{0}$.
The maximum number of iterations was set to 1000.

The numerical experiments were executed on a Dell Precision T7400
(Intel Xeon E5420, 2.5 GHz CPU, 16 GB RAM) running
the Cent OS (kernel 2.6.18)
and the Intel icc 10.1, ifort 10.1 compiler.

In all tests,
ILU(0) was adopted as the preconditioning operation for each of the PCGS algorithms;
here, the value ``zero'' means the {\it fill-in} level.
The ISRVs were set as
$\vc{r}^\flat_0 = \vc{r}_0$ in the conventional PCGS (\Alg~\ref{alg:cgs_precb}),
$\vc{r}^\SH_0 = \vc{r}^+_0$ in the left-PCGS (\Alg~\ref{alg:cgs_Llprec}),
and
$\vc{r}^\SH_0 = \Pinv\vc{r}_0$ in Improved1 and Improved2 
(Algorithms \ref{alg:improved1_cgs_precb} and
\ref{alg:improved2_cgs_precb}, respectively)%
\footnote{
Improved2 was implemented
as the conventional PCGS
with $\vc{r}^\flat_0 = \Pinvt\Pinv \vc{r}_0$ (ISRV9).}
.

We considered the following three cases:
\begin{description}
\item[{\rm (a)}] 
Evaluating the
algorithm relative residual
(see \Fig~\ref{fig:sherman4-ARR}, \ref{fig:watt__1-ARR},
and \Tab\ref{tab:num_results-ARR});

\item[{\rm (b)}] 
Evaluating the
true relative residual
(see \Fig~\ref{fig:sherman4-TRR}, \ref{fig:watt__1-TRR}, and
\Tab\ref{tab:num_results-TRR});
\end{description}
when we have prior knowledge of the exact solution ($\vc{x}_{\rm exact}$);
\begin{description}
\item[{\rm (c)}] 
Evaluating the
true relative error
(see \Fig~\ref{fig:sherman4-TRE}, \ref{fig:watt__1-TRE}, and
\Tab\ref{tab:num_results-TRE}).
\end{description}

We adopted the following stopping criteria:
For case (a),
we adopted the 2-norm of (\ref{eqn:conv_judge}) for 
Algorithms~\ref{alg:cgs_precb}, \ref{alg:improved1_cgs_precb}, and
\ref{alg:improved2_cgs_precb}, and we adopted the
2-norm of (\ref{eqn:conv_judge2}) for
Algorithm~\ref{alg:cgs_Llprec}.
For case (b),
we adopted 
$\|\vc{b}-A\vc{x}_{k+1}\|_2/||\vc{b}||_2 \leq \varepsilon$
for all algorithms.
For case (c),
we adopted 
$\|\vc{x}_{k+1}-\vc{x}_{\rm exact}\|_2/||\vc{x}_{\rm exact}||_2 \leq \varepsilon$
for all algorithms.
We set $\varepsilon=10^{-12}$ for all cases.

We will first focus on the results of the conventional PCGS (\Alg~\ref{alg:cgs_precb}),
as shown in Tables \ref{tab:num_results-ARR} to \ref{tab:num_results-TRE}.
Breakdown occurs for {\tt jpwh\UB 991},
and stagnation occurs for {\tt olm5000} at pitifully insufficient accuracy%
\footnote{
The row marked {\tt olm5000} in all tables contains the results
after 1000 iterations;
furthermore, {\tt olm5000} by \Alg~\ref{alg:cgs_precb} stagnated
after 5000 iterations, due to the size of the matrix.
}%
,
although the other three algorithms
(Algorithms \ref{alg:cgs_Llprec} to \ref{alg:improved2_cgs_precb}) 
were able to solve them.
Note that the disadvantages of the conventional PCGS have already been shown
in \cite{itoh2015a}, and \App~\ref{appsec:sesna_pcgs}.

Next,
we focus on the results of the left-PCGS (\Alg~\ref{alg:cgs_Llprec}),
as shown in Table \ref{tab:num_results-ARR}.
The accuracies of {\tt olm5000} and {\tt watt\UB \UB 1} were highest
for the four PCGS algorithms,
but this has no theoretical underpinning,
because these high accuracies occur %fortuitously
by using the stopping criterion (\ref{eqn:conv_judge2}):
$\|\vc{r}^+_{k+1}\|_2 \; / \; \|\Pinv\vc{b}\|_2$.
In contrast,
the accuracies of {\tt viscoplastic2} was lower
than the two improved PCGS algorithms,
too early convergence occurred.
In many cases,
the left-PCGS algorithm causes the problem of superficial convergence;
see \App~\ref{appsec:sesna_pcgs}.

%%%%%%%%%%%%%%%%%%%%%%%%%%%%%%%%%%%%%%%%%%%%%%%%%%%%%%%%%%%%%%%%%%%%%%%%%
\begin{table}[t]
\caption{(a) Numerical evaluation using the relative residual of each algorithm.
N is the problem size,
and
NNZ is the number of nonzero elements.
The three numbers in each row for the column for each method are as follows: 
the leftmost number is the true relative residual $log_{10}$ 2-norm,
the number in parentheses is the number of iterations required to reach convergence,
and the lower number is the true relative error $log_{10}$ 2-norm.
}
\begin{center}
\begin{small}
\scalebox{0.9}[1.0]{
\begin{tabular}{l|r|r|c|c|c|c}
\hline
\hline
\multicolumn{1}{c|}{\raisebox{1.5ex}{Matrix}}&
\multicolumn{1}{c|}{\raisebox{1.5ex}{N}}&
\multicolumn{1}{c|}{\raisebox{1.5ex}{NNZ}}&
\shortstack{
\\Conventional\\
(\Alg~\ref{alg:cgs_precb})
}
&
\shortstack{
\\Left-PCGS\\
(\Alg~\ref{alg:cgs_Llprec})
}
&
\shortstack{
\\Improved1\\
(\Alg~\ref{alg:improved1_cgs_precb})
}
&
\shortstack{
\\Improved2\\
(\Alg~\ref{alg:improved2_cgs_precb})
}
\\
\hline
add32&4960&19848
&\shortstack{\\-12.17 (35) \\-12.17}
&\shortstack{\\-13.06 (37) \\-12.96}
&\shortstack{\\-12.04 (35) \\-11.96}
&\shortstack{\\-12.04 (35) \\-11.96}
\\
\hline
bfw782a&782&7514
&\shortstack{\\-9.34 (93) \\-10.44}
&\shortstack{\\-12.19 (84) \\-12.23}
&\shortstack{\\-12.19 (84) \\-12.22}
&\shortstack{\\-12.08 (75) \\-11.71}
\\
\hline
jpwh\UB 991&991&6027
&\shortstack{\raisebox{1.5ex}{Breakdown}}
&\shortstack{\\-11.83 (15) \\-12.10}
&\shortstack{\\-12.44 (16) \\-12.53}
&\shortstack{\\-12.44 (16) \\-12.53}
\\
\hline
olm5000&5000&19996
&\shortstack{\\-0.18 (Stag.) \\4.22}
&\shortstack{\\-12.79 (38) \\-10.64}
&\shortstack{\\-12.20 (34) \\-8.05}
&\shortstack{\\-12.21 (33) \\-8.00}
\\
\hline
poisson3Db&85623&2374949
&\shortstack{\\-10.14 (122) \\-10.33}
&\shortstack{\\-12.93 (119) \\-13.31}
&\shortstack{\\-12.49 (123) \\-13.39}
&\shortstack{\\-11.79 (117) \\-12.07}
\\
\hline
sherman4&1104&3786
&\shortstack{\\-12.69 (34) \\-13.83}
&\shortstack{\\-11.68 (32) \\-12.82}
&\shortstack{\\-12.69 (33) \\-13.82}
&\shortstack{\\-12.69 (33) \\-13.83}
\\
\hline
viscoplastic2&32769&381326
&\shortstack{\\-7.55 (812) \\-4.68}
&\shortstack{\\-10.26 (775) \\-7.54}
&\shortstack{\\-11.80 (844) \\-8.69}
&\shortstack{\\-11.81 (886) \\-8.84}
\\
\hline
watt\UB \UB 1&1856&11360
&\shortstack{\\-13.01 (27) \\-5.96}
&\shortstack{\\-15.48 (41) \\-12.63}
&\shortstack{\\-12.11 (35) \\-9.77}
&\shortstack{\\-12.11 (35) \\-9.77}
\\
\hline
\end{tabular}
\label{tab:num_results-ARR}
}
\end{small}
\vspace{-15pt}
\end{center}
\end{table}

%%%%%%%%%%%%%%%%%%%%%%%%%%%%%%%%%%%%%%%%%%%%%%%%%%%%%%%%%%%%%%%%%%%%%%%%%

%%%%%%%%%%%%%%%%%%%%%%%%%%%%%%%%%%%%%%%%%%%%%%%%%%%%%%%%%%%%%%%%%%%%%%%%%
\begin{table}[thp]
\caption{(b) Numerical evaluation using the true relative residual of each algorithm.}
\begin{center}
\begin{small}
\scalebox{0.9}[1.0]{
\begin{tabular}{l|r|r|c|c|c|c}
\hline
\hline
\multicolumn{1}{c|}{\raisebox{1.5ex}{Matrix}}&
\multicolumn{1}{c|}{\raisebox{1.5ex}{N}}&
\multicolumn{1}{c|}{\raisebox{1.5ex}{NNZ}}&
\shortstack{
\\Conventional\\
(\Alg~\ref{alg:cgs_precb})
}
&
\shortstack{
\\Left-PCGS\\
(\Alg~\ref{alg:cgs_Llprec})
}
&
\shortstack{
\\Improved1\\
(\Alg~\ref{alg:improved1_cgs_precb})
}
&
\shortstack{
\\Improved2\\
(\Alg~\ref{alg:improved2_cgs_precb})
}
\\
\hline
add32
&4960&19848
&\shortstack{\\-12.17 (35) \\-12.17}
&\shortstack{\\-12.04 (35) \\-11.96}
&\shortstack{\\-12.04 (35) \\-11.96}
&\shortstack{\\-12.04 (35) \\-11.96}
\\
\hline
bfw782a&782&7514
&\shortstack{\\-9.34 (Stag.) \\-10.44}
&\shortstack{\\-12.19 (84) \\-12.23}
&\shortstack{\\-12.19 (84) \\-12.22}
&\shortstack{\\-12.08 (75) \\-11.71}
\\
\hline
jpwh\UB 991
&991&6027
&\shortstack{\raisebox{1.5ex}{Breakdown}}
&\shortstack{\\-12.44 (16) \\-12.53}
&\shortstack{\\-12.44 (16) \\-12.53}
&\shortstack{\\-12.44 (16) \\-12.53}
\\
\hline
olm5000
&5000&19996
&\shortstack{\\-0.18 (Stag.) \\4.22}
&\shortstack{\\-12.49 (31) \\-8.28}
&\shortstack{\\-12.20 (34) \\-8.05}
&\shortstack{\\-12.21 (33) \\-8.00}
\\
\hline
poisson3Db
&85623&2374949
&\shortstack{\\-10.14 (Stag.) \\-10.33}
&\shortstack{\\-12.08 (113)   \\-12.95}
&\shortstack{\\-12.49 (123)   \\-13.39}
&\shortstack{\\-11.77 (Stag.) \\-12.06}
\\
\hline
sherman4
&1104&3786
&\shortstack{\\-12.69 (34) \\-13.83}
&\shortstack{\\-12.68 (33) \\-13.81}
&\shortstack{\\-12.69 (33) \\-13.82}
&\shortstack{\\-12.69 (33) \\-13.83}
\\
\hline
viscoplastic2&32769&381326
&\shortstack{\\-7.55 (Stag.) \\-4.68}
&\shortstack{\\-10.27 (Stag.) \\-8.01}
&\shortstack{\\-11.84 (Stag.) \\-8.95}
&\shortstack{\\-11.82 (Stag.) \\-8.93}
\\
\hline
watt\UB \UB 1
&1856&11360
&\shortstack{\\-13.01 (27) \\-5.96}
&\shortstack{\\-12.11 (35) \\-9.77}
&\shortstack{\\-12.11 (35) \\-9.77}
&\shortstack{\\-12.11 (35) \\-9.77}
\\
\hline
\end{tabular}
\label{tab:num_results-TRR}
}
\end{small}
\vspace{10pt}
%%%%%%%%%%%%%%%%%%%%%%%%%%%%%%%%%%%%%%%%%%%%%%%%%%%%%%%%%%%%%%%%%%%%%%%%%
%
%
%%%%%%%%%%%%%%%%%%%%%%%%%%%%%%%%%%%%%%%%%%%%%%%%%%%%%%%%%%%%%%%%%%%%%%%%%
\caption{(c) Numerical evaluation using the true relative error of each algorithm.}
\vspace{7pt}
\begin{small}
\scalebox{0.9}[1.0]{
\begin{tabular}{l|r|r|c|c|c|c}
\hline
\hline
\multicolumn{1}{c|}{\raisebox{1.5ex}{Matrix}}&
\multicolumn{1}{c|}{\raisebox{1.5ex}{N}}&
\multicolumn{1}{c|}{\raisebox{1.5ex}{NNZ}}&
\shortstack{
\\Conventional\\
(\Alg~\ref{alg:cgs_precb})
}
&
\shortstack{
\\Left-PCGS\\
(\Alg~\ref{alg:cgs_Llprec})
}
&
\shortstack{
\\Improved1\\
(\Alg~\ref{alg:improved1_cgs_precb})
}
&
\shortstack{
\\Improved2\\
(\Alg~\ref{alg:improved2_cgs_precb})
}
\\
\hline
add32
&4960&19848
&\shortstack{\\-12.17 (35) \\-12.17}
&\shortstack{\\-12.00 (36) \\-12.29}
&\shortstack{\\-12.00 (36) \\-12.29}
&\shortstack{\\-12.00 (36) \\-12.29}
\\
\hline
bfw782a&782&7514
&\shortstack{\\-9.34 (Stag.) \\-10.44}
&\shortstack{\\-12.19 (84) \\-12.23}
&\shortstack{\\-12.19 (84) \\-12.22}
&\shortstack{\\-12.20 (84) \\-12.23}
\\
\hline
jpwh\UB 991
&991&6027
&\shortstack{\raisebox{1.5ex}{Breakdown}}
&\shortstack{\\-11.83 (15) \\-12.10}
&\shortstack{\\-11.83 (15) \\-12.10}
&\shortstack{\\-11.83 (15) \\-12.10}
\\
\hline
olm5000
&5000&19996
&\shortstack{\\-0.18 (Stag.) \\4.22}
&\shortstack{\\-12.79 (Stag.) \\-11.23}
&\shortstack{\\-12.80 (49) \\-13.22}
&\shortstack{\\-12.59 (52) \\-13.09}
\\
\hline
poisson3Db
&85623&2374949
&\shortstack{\\-10.14 (Stag.) \\-10.33}
&\shortstack{\\-11.24 (111)   \\-12.21}
&\shortstack{\\-11.61 (117)   \\-12.57}
&\shortstack{\\-11.53 (116)   \\-12.04}
\\
\hline
sherman4
&1104&3786
&\shortstack{\\-12.69 (34) \\-13.83}
&\shortstack{\\-11.68 (32) \\-12.82}
&\shortstack{\\-11.68 (32) \\-12.82}
&\shortstack{\\-11.68 (32) \\-12.82}
\\
\hline
viscoplastic2&32769&381326
&\shortstack{\\-7.55 (Stag.) \\-4.68}
&\shortstack{\\-10.27 (Stag.) \\-8.01}
&\shortstack{\\-11.84 (Stag.) \\-8.95}
&\shortstack{\\-11.82 (Stag.) \\-8.93}
\\
\hline
watt\UB \UB 1
&1856&11360
&\shortstack{\\-18.06 (41) \\-12.13}
&\shortstack{\\-14.28 (40) \\-12.02}
&\shortstack{\\-14.26 (40) \\-12.01}
&\shortstack{\\-14.26 (40) \\-12.05}
\\
\hline
\end{tabular}
\label{tab:num_results-TRE}
}
\end{small}
\end{center}
\end{table}
%%%%%%%%%%%%%%%%%%%%%%%%%%%%%%%%%%%%%%%%%%%%%%%%%%%%%%%%%%%%%%%%%%%%%%%%%

Next,
it is very important to compare
cases (a) and (b)  (Tables \ref{tab:num_results-ARR} and \ref{tab:num_results-TRR})
with
case (c) (\Tab\ref{tab:num_results-TRE}),
in order to determine the crucial ways in which they differ.
Because (a) and (b)
can be evaluated without knowing the exact solution
but
(c) requires the exact solution,
it is important to examine the results when the exact solution is known.
Comparing the results
for {\tt bfw782a}, {\tt poisson3Db}, {\tt viscoplastic2}, and {\tt watt\UB\UB 1}
in cases (a) and (b) (Tables \ref{tab:num_results-ARR} and \ref{tab:num_results-TRR}),
the conventional PCGS (\Alg~\ref{alg:cgs_precb})
has results in which the true relative residual or true relative error
(or both)
is much less accurate than those obtained by the other algorithms,
and only in the conventional PCGS does stagnation occur
at insufficient accuracy%.
\footnote{
The results of
{\tt poisson3Db} with Improved2 in \Tab\ref{tab:num_results-TRR}
and
{\tt olm5000} with the left-PCGS in \Tab\ref{tab:num_results-TRE}
can be considered to be sufficiently accurate,
because
they had nearly converged to within $10^{-12}$.
We note that $\varepsilon=10^{-12}$ is a stringent value
for the tolerance for the true relative residual and the true relative error.
}%
.
In particular,
the conventional PCGS is the fastest to converge
 for {\tt watt\UB\UB 1} in
cases (a) and (b)
(Tables \ref{tab:num_results-ARR} and \ref{tab:num_results-TRR}),
but
this is undesirable,
because
when convergence occurs too quickly,
the relative residual and the true relative residual may fail to meet the criterion for accuracy.
On the other hand, based on
the true relative error
in case (c) (\Tab\ref{tab:num_results-TRE}), we see that
the conventional PCGS converges after almost the same number of iterations
as do the other methods.

%%%%%%%%%%%%%%%%%%%%%%%%%%%%%%%%%%%%%%%%%%%%%%%%%%%%%%%%%%%%%%%%%%%%%%%%%
\def\CASE{sherman4}
\begin{figure}[!tp]
\begin{center}
\resizebox*{0.6\columnwidth}{0.6\height}{
\includegraphics*{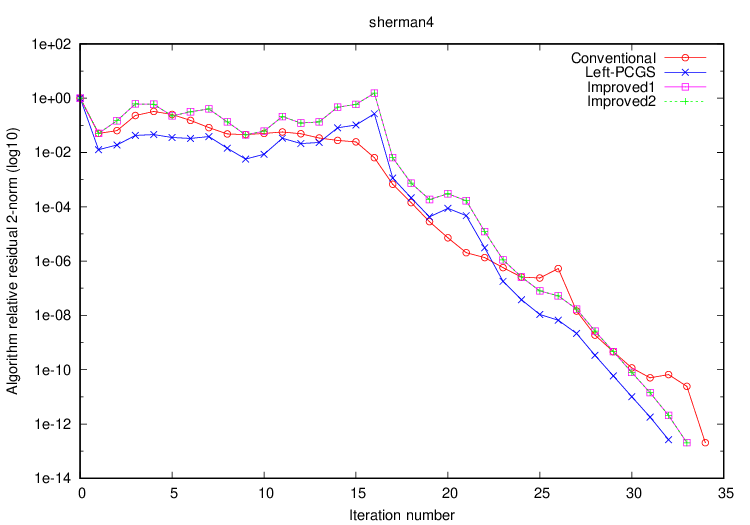}}
\caption{(a) Convergence history of the algorithm relative residual 2-norm for each of the four algorithms (\CASE).}
\label{fig:sherman4-ARR}
\end{center}
\def\CASE{sherman4}
\begin{center}
\resizebox*{0.6\columnwidth}{0.6\height}{
\includegraphics*{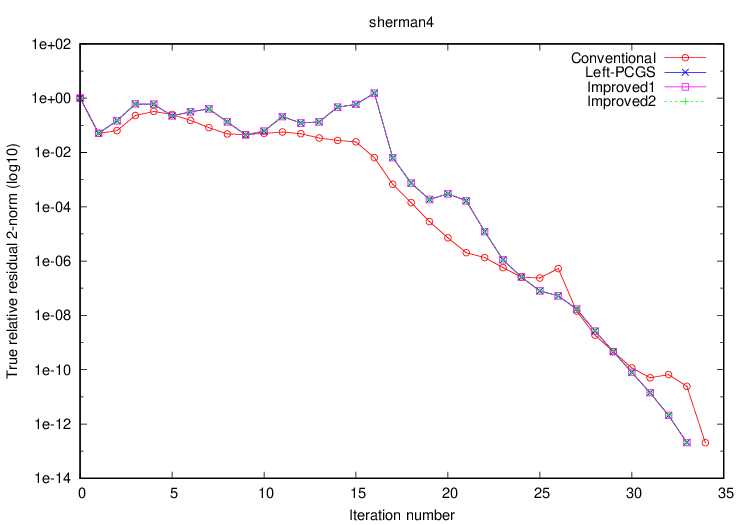}}
\caption{(b) Convergence history of the true relative residual 2-norms for each of the four algorithms (\CASE).}
\label{fig:sherman4-TRR}
\end{center}
\def\CASE{sherman4}
\begin{center}
\resizebox*{0.6\columnwidth}{0.6\height}{
\includegraphics*{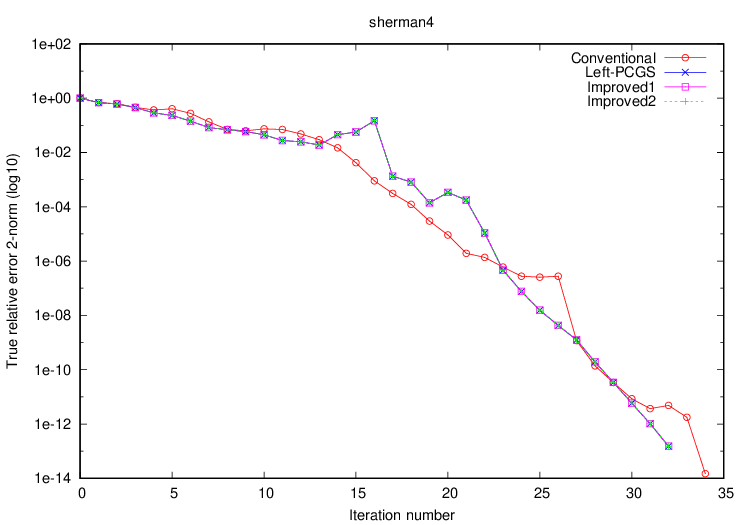}}\caption{(c) Convergence history of the true relative error 2-norm for each of the four algorithms (\CASE).}
\label{fig:sherman4-TRE}
\end{center}
\end{figure}
%
%%%%%%%%%%%%%%%%%%%%%%%%%%%%%%%%%%%%%%%%%%%%%%%%%%%%%%%%%%%%%%%%%%%%%%%%%

%%%%%%%%%%%%%%%%%%%%%%%%%%%%%%%%%%%%%%%%%%%%%%%%%%%%%%%%%%%%%%%%%%%%%%%%%
\def\CASE{watt\UB \UB 1}
\begin{figure}[!tp]
\begin{center}
\resizebox*{0.6\columnwidth}{0.6\height}{
\includegraphics*{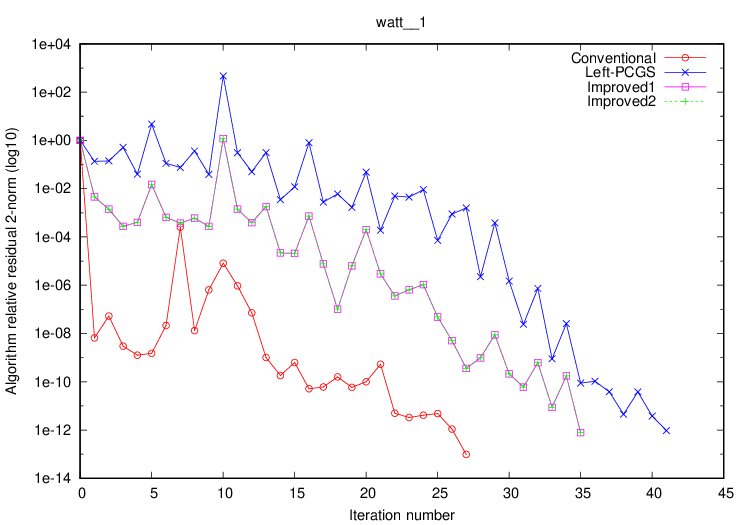}}
\caption{(a) Convergence history of the algorithm relative residual 2-norm for each of the four algorithms (\CASE).}
\label{fig:watt__1-ARR}
\end{center}
\def\CASE{watt\UB \UB 1}
\begin{center}
\resizebox*{0.6\columnwidth}{0.6\height}{
\includegraphics*{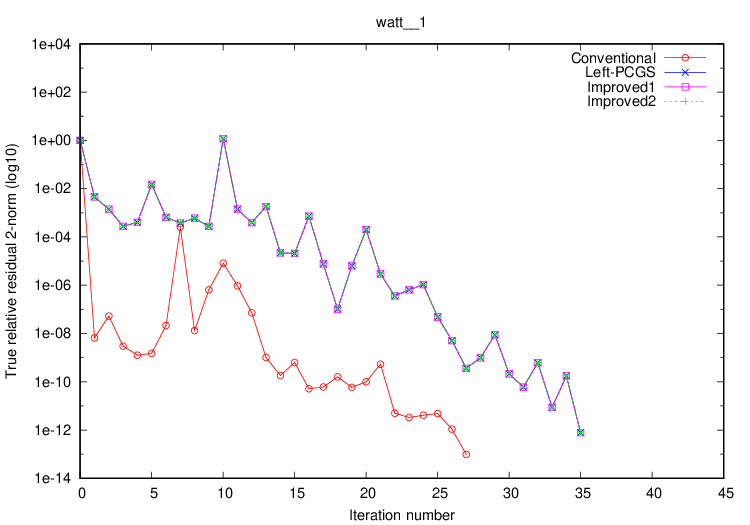}}
\caption{(b) Convergence history of the true relative residual 2-norm for each of the four algorithms (\CASE).}
\label{fig:watt__1-TRR}
\end{center}
\def\CASE{watt\UB \UB 1}
\begin{center}
\resizebox*{0.6\columnwidth}{0.6\height}{
\includegraphics*{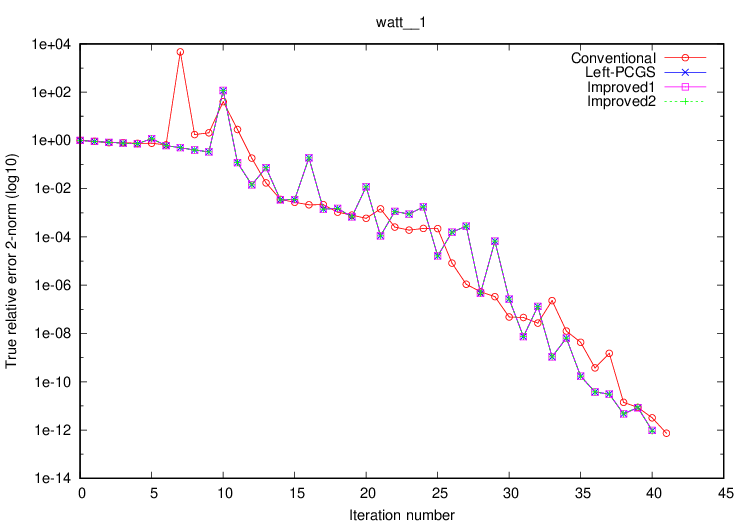}}
\caption{(c) Convergence history of the true relative error 2-norm for each of the four algorithms (\CASE).}
\label{fig:watt__1-TRE}
\end{center}
\end{figure}
%
%%%%%%%%%%%%%%%%%%%%%%%%%%%%%%%%%%%%%%%%%%%%%%%%%%%%%%%%%%%%%%%%%%%%%%%%%

From the graphs in
Figures \ref{fig:sherman4-ARR} to \ref{fig:watt__1-TRE},
we can see the following:
in case (a), Improved1, Improved2, and the left-PCGS show different convergence behaviors, but in cases (b) and (c), they show similar behaviors.
These results correspond to the analysis
in Section~\ref{ssec:analysis_4pcgs}.
Therefore,
Algorithms \ref{alg:improved1_cgs_precb} and \ref{alg:improved2_cgs_precb}
are coordinative to \Alg~\ref{alg:cgs_Llprec} regarding the structures of
the solution vector for the polynomial ${\it\Phi}_k(\tA)\tvc{r}_0$,
despite the difference between the residual vectors
$\vc{r}^+_k$ for the left-PCGS (\Alg~\ref{alg:cgs_Llprec})
and
$\vc{r}_k$ for Improved1 and Improved2 
 (Algorithms \ref{alg:improved1_cgs_precb} and \ref{alg:improved2_cgs_precb},
 respectively).
The conventional PCGS had a convergence behavior
that differs from those of all of the other algorithms for cases (a) to (c).

These numerical results conform to the behavior expected based on
the discussion of the relation between
the structure of the solution vector and the polynomial.
We compared the numerical results with the theoretical results of
Sections~\ref{sssec:comparison_conventional_pcgs} and~\ref{ssec:analysis_4pcgs},
and these results are summarized as follows:
\begin{enumerate}
\item
For case (a), the difference between
 the residual vector $\vc{r}^+_k$ of the left-PCGS 
and
 $\vc{r}_k$
has been verified.
Accuracy of the left-PCGS is not conclusive,
because its high accuracy is caused not by necessity
but by accident of its stopping criterion:
$\|\vc{r}^+_{k+1}\|_2 / \|\Pinv\vc{b}\|_2$.

\item
For cases (b) and (c), we verified (\ref{eqn:structure_of_solution}):

$\vc{x}_k = \vc{x}_\km + {\it\Phi}^{\rm L}_\km(\Pinv A)\vc{r}^+_0
% $\vc{x}_k = \vc{x}_\km + \alpha^{\rm L}_\km{\it\Phi}^{\rm L}_\km(\Pinv A)\vc{r}^+_0
\SP\SP \mapsto \SP\SP
 \vc{x}_k = \vc{x}_\km + \Pinv{\it\Phi}^{\rm L}_\km(A\Pinv)\vc{r}_0$.
%  \vc{x}_k = \vc{x}_\km + \alpha^{\rm L}_\km\Pinv{\it\Phi}^{\rm L}_\km(A\Pinv)\vc{r}_0$.

Furthermore,
we verified (\ref{eqn:differ_bothKSP}):
${\it\Phi}^{\rm R}_\km(A\Pinv)\vc{r}_0
 \neq
 {\it\Phi}^{\rm L}_\km(A\Pinv)\vc{r}_0$,

\item
The differences between the conventional PCGS, the left-PCGS,
Improved1, and Improved2 have been confirmed by examining their convergence behavior.
In other word, we considered the relation of the solution vector
and the polynomial ${\it\Phi}_k(\tA)\tvc{r}_0$
between
 the right system (the conventional PCGS)
and
 the left-PCGS, and between the coordinative PCGSs (Improved1 and Improved2) and the left-PCGS.
\end{enumerate}

\subsection{Behavior of the PCGS when it is switched by the ISRV}
\label{ssec:num_exp2}

In this subsection,
the experimental environment was the same as that described in Section~\ref{ssec:num_exp1},
except that we used Matlab 7.8.0 (R2009a),
and
we gave different ISRVs to the conventional PCGS and Improved1.

We compared five different PCGS algorithms, including using a different ISRV.
In the figures, we use the following labels:
``Conventional'' means the conventional PCGS
(\Alg~\ref{alg:cgs_precb}), for which the ISRV is $\vc{r}^\flat_0=\vc{r}_0$;
 this is a right-preconditioned system.
``Impr1-ISRV1'' means Improved1 (\Alg~\ref{alg:improved1_cgs_precb}) with ISRV1:
$\vc{r}^\SH_0=\Pinv\vc{r}_0$.
``Impr1-ISRV2'' means Improved1 with ISRV2: $\vc{r}^\SH_0=\PO^\T\vc{r}_0$.
``Left'' means
 the left-PCGS (\Alg~\ref{alg:cgs_Llprec}),
 for which the ISRV is $\vc{r}^\SH_0=\vc{r}^+_0$.
``Conv\UB{}ISRV9'' means the conventional PCGS with ISRV9:
 $\vc{r}^\flat_0=\Pinvt\Pinv\vc{r}_0$.

\def\FIGSIZE{0.7}
\def\RATIOTWO{0.7}
\def\CASE{sherman4}
\begin{figure}[thp]
\begin{center}
\resizebox*{\FIGSIZE\columnwidth}{\RATIOTWO\height}{
\includegraphics*{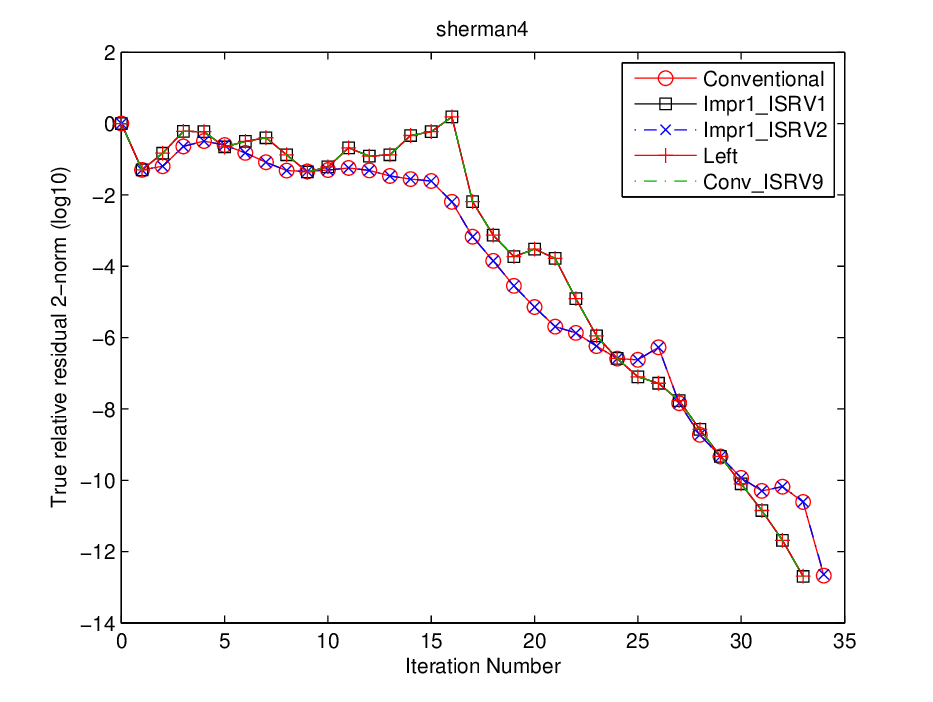}}
\caption{ Convergence history of the true relative residual 2-norm
 of the right- and left-preconditioned PCGS,
 for each of the five PCGS algorithms with ISRV switching (\CASE).}
\label{fig:sherman4-isrv}
\end{center}
\def\CASE{watt\UB \UB 1}
\begin{center}
\resizebox*{\FIGSIZE\columnwidth}{\RATIOTWO\height}{
\includegraphics*{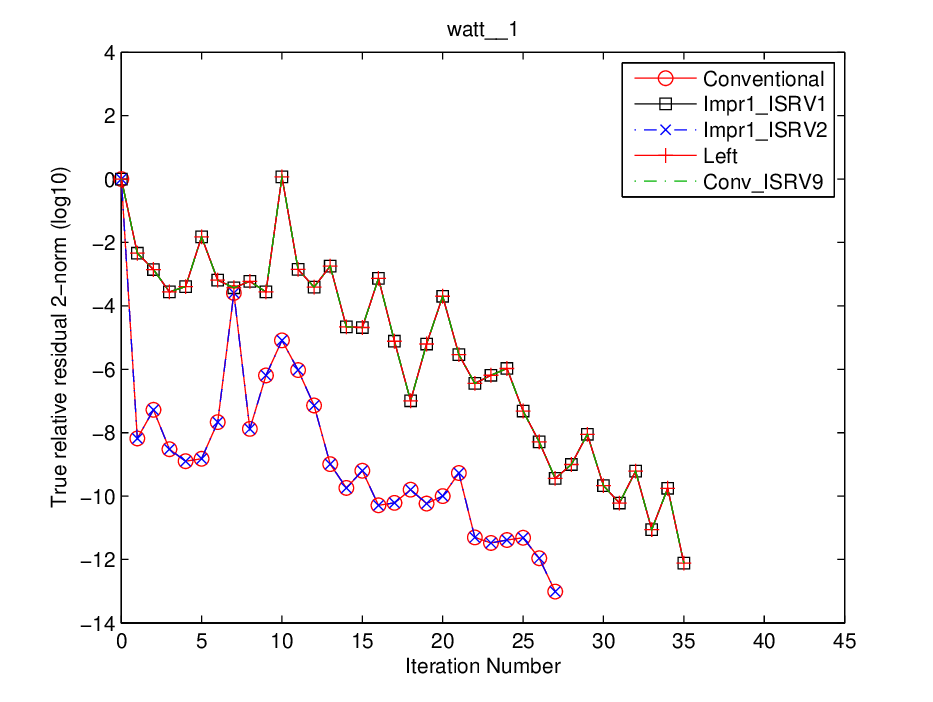}}
\caption{Convergence history of true relative residual 2-norm
 of the right- and left-preconditioned PCGS,
 for each of the five PCGS algorithms with ISRV switching (\CASE).}
\label{fig:watt__1-isrv}
\end{center}
\end{figure}

The convergence histories of
``Conventional,'' ``Impr1-ISRV1,'' and ``Left''
 in Figures \ref{fig:sherman4-isrv} and \ref{fig:watt__1-isrv}
 are the same as those of
``Conventional,'' ``Improved1,'' and ``Left-PCGS,'' respectively, 
in Figures \ref{fig:sherman4-TRR} and \ref{fig:watt__1-TRR}.

In both figures,
``Impr1-ISRV2'' and ``Conv\UB{}ISRV9'' were added
to verify \Thm~\ref{thm:isrv_switch}.
The convergence history of ``Impr1-ISRV2'' is the same as that of ``Conventional,''
and those of ``Impr1-ISRV1'' and ``Conv\UB{}ISRV9'' are the same as that of ``Left.''

We have numerically verified the claims of Section~\ref{sec:congruent_and_direction_pcgs};
in particular, we have verified
\Thm~\ref{thm:isrv_switch}.

%%%%%%  section  %%%%%%%%%%%%%%%%%%%%%%%%%%%%%%%%%%%%%%%%%%%%%%%%%%%%%%%%%%%%%%%%%
\section{Conclusions}
\label{sec:conclusion}

In this paper,
an improved PCGS algorithm \cite{itoh2015a} has been analyzed
by mathematically comparing four different PCGS algorithms,
and we have focused on the structures of the solution vector
 and
their polynomial ${\it\Phi}_k(\tA)\tvc{r}_0$.
From our analysis and numerical results,
we have verified two improved PCGS algorithms.
They are both coordinative to the left-preconditioned systems,
although
their residual vector maintains the basic form $\vc{r}_k$,
not $\vc{r}^+_k$.
For both algorithms,
the structures of the solution vector and the polynomial are
$\vc{x}_k = \vc{x}_\km + \Pinv{\it\Phi}^{\rm L}_\km(A\Pinv)\vc{r}_0$.
Furthermore,
the numerical results of
the improved PCGS with the ILU(0) preconditioner
show many advantages,
such as effectiveness and consistency across several preconditioners,
have also been shown; see \cite{itoh2015a} and \App~\ref{appsec:sesna_pcgs}.
We note that the improved PCGS algorithms share some of the advantages
of the conventional PCGS (the right-preconditioned system)
and the left-PCGS algorithms,
while they avoid some of their disadvantages.
Accuracy of the left-PCGS is not controllable,
but the Improved1 PCGS algorithm has availability of
further improving on accuracy by controlling both residual vectors
$\vc{r}_{k+1}$ and $\Pinv\vc{r}_{k+1} \; ( \; =\vc{r}^+_{k+1} \; )$
in algorithm.
This is our future assignment.

We presented a general definition
of the direction of a preconditioned system of linear equations.
Furthermore,
we have shown that the direction of a preconditioned system for CGS
is switched by the construction and setting of
the ISRV.
This is because the direction of the preconditioning conversion
is congruent.
We have also shown that
the direction of a preconditioned system for CGS is determined
by the operations of $\alpha_k$ and $\beta_k$
and 
that these intrinsic operations are based on biorthogonality
and biconjugacy.
However,
the structures of these intrinsic operations
are the same in all four of the PCGS algorithms.
Therefore,
we have focused on the ability of the ISRV to
switch the direction of a preconditioned system,
and such a mechanism may be unique to the bi-Lanczos-type algorithms
that are based on the BiCG method;
for example,
preconditioned BiCGStab, preconditioned GPBiCG and so on.
Here,
there is an arguable issue on designing their algorithm constructions
because such preconditioned bi-Lanczos-type methods
have minimal residual operators
that have no congruency of preconditioning conversion.
This is also our future assignment.

As we analyzed the four PCGS algorithms, we paid particular attention to the
vectors.
We note that 
there exist preconditioned BiCG (PBiCG) algorithms that
correspond to the preconditioning conversion of each of the PCGS algorithms.
The polynomial structure of the PBiCG can be minutely analyzed 
by replacing the vectors of the PCGS.
We have analyzed the four PBiCG algorithms in parallel \cite{itoh2016b},
and each PBiCG corresponds to one of the four PCGS algorithms
in this paper.
In \cite{itoh2016b}, using the ISRV to
switch the direction of a preconditioned system was
discussed in detail.

%%%%%%%%%%%%%%%%%%%%%%%%%%%%%%%%%%%%%%%%%%%%%%%%%%%%%%%%%%%%%%%%%%%%%%%%%%%%%%%%%%

%%%% Acknowledgments %%%%%%%%
\section*{Acknowledgments}
This work was partially supported by JSPS KAKENHI Grant Number JP25390145, JP18K11342.

%%%% Bibliography  %%%%%%%%%%

%%%%%%%%%%%%%%%%%%%%%%%%%%%%%%%%%%%%%%%%%%%%%%%%%%%%%%%%%%%%%%%%%%%%%%%%%%%%%%%%%%
\appendix
\section{Systematic performance evaluation of PCGS}
\label{appsec:sesna_pcgs}

\begin{figure}[htbp]
\begin{center}
\resizebox*{0.68\columnwidth}{!}{
\includegraphics*{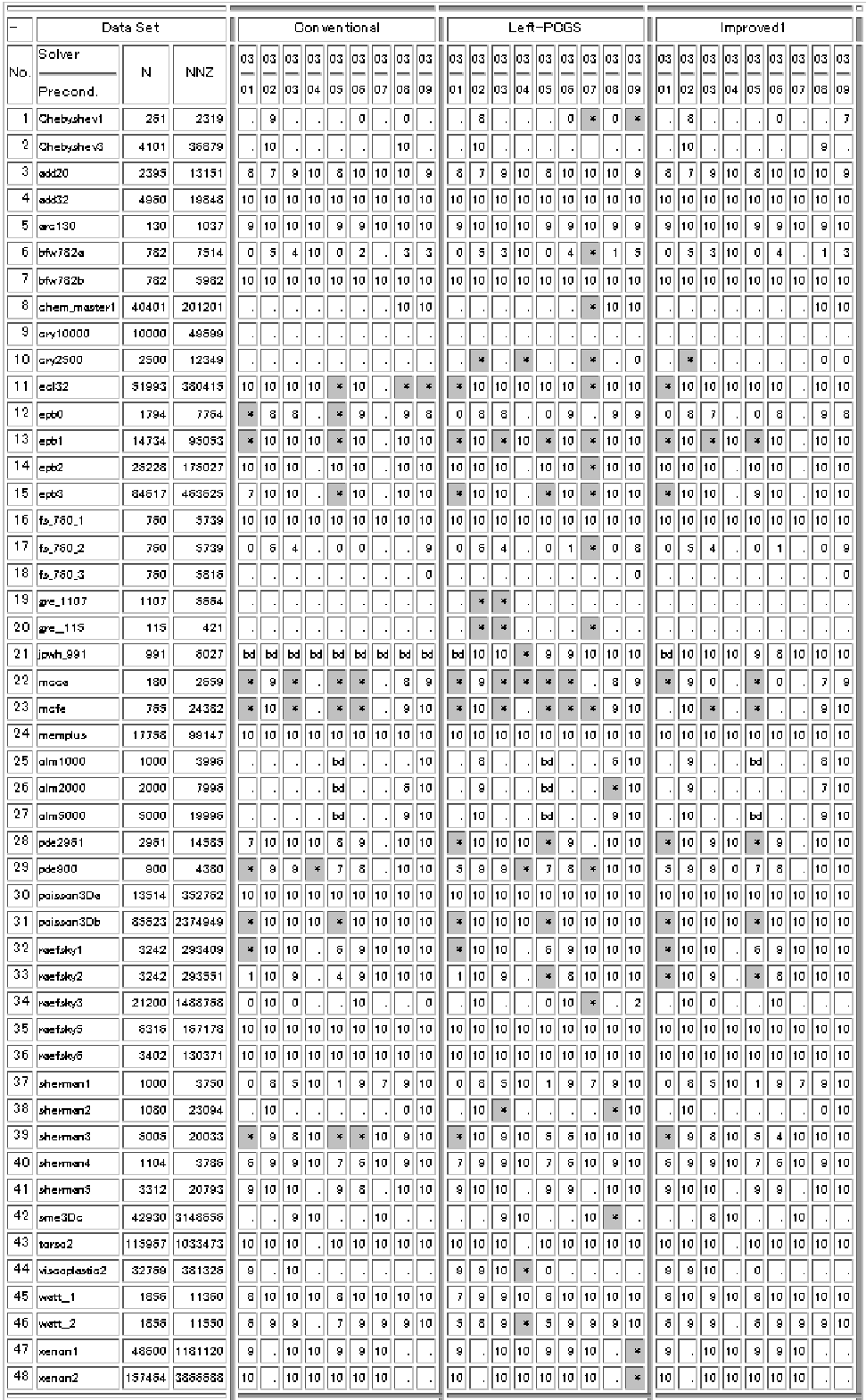}}
\caption{Superficial convergence with
three different PCGS methods and nine different preconditioners.
Superficial convergence is indicated by an asterisk {\tt *} with a gray background.
Thick lines group columns for (from left to right) the
conventional PCGS (\Alg~\ref{alg:cgs_precb}), the
left-PCGS (\Alg~\ref{alg:cgs_Llprec}), and
Improved1 (\Alg~\ref{alg:improved1_cgs_precb}). Each such
frame has a cell for each of the nine preconditioners (from left to right):
Point Jacobi {\tt [01]}, ILU(0) {\tt [02]}, SSOR {\tt [03]},
Hybrid type {\tt [04]}, I+S {\tt [05]}, SAINV {\tt [06]},
SA-AMG {\tt [07]} (with ``{\tt -saamg\UB{}unsym true}'' option),
Crout ILU {\tt [08]}, and ILUT {\tt [09]},
where the two digits in {\tt [ ]} indicate the identification number of the solver or preconditioner.
{\tt N} is the matrix size, and {\tt NNZ} is the number of nonzero elements.}
\label{fig:ConvLeftImpr1_woList}
\end{center}
\end{figure}

We consider the superficial convergence
on three types of preconditioned CGS algorithms;
the results shown in \Fig~\ref{fig:ConvLeftImpr1_woList} provide an overview
of this weakness in the left-PCGS.
The figure presents
a systematic performance evaluation of solving a linear equation
by three PCGS algorithms
(conventional, left, and improved1)
and with nine kinds of preconditioners.
Here,
we briefly discuss the figure;
the details of a systematic performance evaluation can be found
in \Ref~\cite{itoh2010e},
which discusses the right- and the left-preconditioned systems.

The conditions of the numerical experiments presented in this appendix
were almost the same as those in Section~\ref{ssec:num_exp1}
and \Ref~\cite[Section 9.2.2]{itoh2010e},
except for the computational environment.
The systematic performance evaluation was executed on
a Hitachi HA8000 (AMD Opteron 8356; 2.3 GHz CPU; memory size: 32 GB/node)
running the Red Hat Enterprise Linux 5 and
the Intel icc 10.1, ifort 10.1 compiler. The
maximum number of iterations was equal to the size of each matrix.
The test problems were also generated as 
in Section~\ref{ssec:num_exp1} and \Ref~\cite[Section 9.2.2]{itoh2010e}.

\Fig~\ref{fig:ConvLeftImpr1_woList} shows 
the solution performance data. The rows indicate the
forty-eight kinds of linear equations for the coefficient matrix names that are listed in \Tab\ref{tab:mat_sesna}; the columns indicate the three PCGS algorithms, each of which is subdivided into columns for the nine preconditioners.

\begin{table}[thp]
\caption{Real nonsymmetric matrices used
for systematic performance evaluation in \Fig~\ref{fig:ConvLeftImpr1_woList}.}
\begin{center}
\begin{small}
\scalebox{0.9}[1.0]{
\begin{tabular}{l}
\hline
\hline
Chebyshev\{1,3\},
add\{20,32\},
arc130,
bfw782\{a,b\},
chem\UB master1,
cry\{10000,2500\},
\\
ecl32,
epb[0-3],
fs\UB 760\UB [1-3],
gre\UB \{1107,115\},
jpwh\UB 991,
mcca,
mcfe,
memplus,
\\
olm\{1000, 2000, 5000\},
pde\{2961,900\},
poisson3D\{a,b\},
raefsky\{1,2,3,5,6\},
\\
sherman[1-5],
sme3Dc,
torso2,
viscoplastic2,
watt\UB \UB \{1,2\},
xenon\{1,2\}
\\
\hline
\end{tabular}
}
\end{small}
\end{center}
\label{tab:mat_sesna}
\end{table}

The contents of each cell are as follows.
The number in each cell indicates the convergence rate%
\footnote{
We based the maximum number of iterations on the size of the problem,
and converted this to a percentage to obtain the convergence score.
In this study,
we calculate the score for when the number of iterations required for convergence
is less than or equal to $20\%$ of the matrix size,
then
\begin{eqnarray}
\mbox{\tt score} &=&
10 -
\left[
\frac{\mbox{Required number of iterations {\tt (iter)}} - 1}
     {\mbox{Size of the coefficient matrix {\tt (N)}}}
\times {\tt div}
\right], 
\;\;\;\;
\mbox{(in this study, {\tt div = 50}),}
\NON
\end{eqnarray}
otherwise,
{\tt score = 0}.
If {\tt score < 0}, then {\tt score = 0}.
Here, $\left[ \;\; \right]$ indicates Gauss notation.
If {\tt score = 10},
the number of iterations required for convergence is less than or equal to $2\%$ of the matrix size,
and a lower score means slower convergence.
If {\tt score = 0},
the number of iterations required for convergence is greater than $20\%$ of the matrix size \cite{itoh2010e}.
Example: We solve a linear equation of matrix size {\tt N = 782}.
If {\tt iter = 14, 84, 148, 259}, then {\tt score = 10, 5, 1, 0}, respectively.

However, in this appendix,
it is not our main purpose to discuss the value of the score
 but to compare the instances of superficial convergence ({\tt *})
or other problem cases (`period', `{\tt bd}', and `blank').
}.
A period ({\tt .}) indicates that it did not converge until the maximum number of iterations, ``{\tt bd}'' indicates that the process broke down,
an asterisk ($\ast$) with a gray background indicates that the convergence
was superficial convergence,
and a blank indicates that it was not solved for some other reason.

Here,
we consider the number of cases of superficial convergence.
Superficial convergence occurs when the residual vector implies that convergence has occurred, but 
the solution vector is not a sufficiently accurate approximation to the true residual vector.

The stopping criteria were as follows:
\begin{equation}
\reeqno{\ref{eqn:conv_judge}\sq}
\|\vc{r}_{k+1}\|_2 \; / \; \|\vc{b}\|_2 \leq \varepsilon,
  \;\;\;\;
\mbox{(Conventional and Improved1 PCGS)},
\end{equation}
\vspace*{-18pt}
\begin{equation}
\reeqno{\ref{eqn:conv_judge2}\sq}
\|\vc{r}^+_{k+1}\|_2 \; / \; \|\Pinv\vc{b}\|_2 \leq \varepsilon,
  \;\;\;\;
\mbox{(Left-PCGS)}.
\end{equation}
Here, $\vc{r}_{k+1}$ and $\vc{r}^+_{k+1}$ are the residual vectors
for the corresponding PCGS.
\setcounter{equation}{0}
When these conditions are satisfied,
the numerical solution ($\hat{\vc{x}}_{k+1}$) is obtained,
and
the true relative residual is calculated as follows:
\begin{eqnarray}
&& \|\vc{b} - A \hat{\vc{x}}_{k+1}\|_2 \; / \; \|\vc{b}\|_2 \leq \hat{\varepsilon} .
\label{eqn:conv_judge_TRR_sesna}
\end{eqnarray}
In this systematic performance evaluation,
$\varepsilon$ in 
(\ref{eqn:conv_judge}\sq) and 
(\ref{eqn:conv_judge2}\sq) was set to $1.0 \times 10^{-12}$.
On the other hand,
$\hat{\varepsilon}$ in
(\ref{eqn:conv_judge_TRR_sesna}) was set to $1.0 \times 10^{-8}$,
because $\hat{\varepsilon} = 1.0 \times 10^{-12}$ is
a stringent value for the tolerance for the true relative residual.
We consider that superficial convergence has occurred when 
(\ref{eqn:conv_judge}\sq) or (\ref{eqn:conv_judge2}\sq)
is satisfied, but (\ref{eqn:conv_judge_TRR_sesna})
is not.

%%%%%%%%%%%%%%%%%%%%%%%%%%%%%%%%%%%%%%%%%%%%%%%%%%%%%%%%%%%%%%%%%%%%%%%%
\begin{table}[thp]
\caption{Number of cases of superficial convergence for each PCGS
in \Fig~\ref{fig:ConvLeftImpr1_woList}.}
\begin{center}
\begin{tabular}{c|c|c}
\hline
\hline
\shortstack{
\\Conventional\\
(\Alg~\ref{alg:cgs_precb})
}
&
\shortstack{
\\Left-PCGS\\
(\Alg~\ref{alg:cgs_Llprec})
}
&
\shortstack{
\\Improved1\\
(\Alg~\ref{alg:improved1_cgs_precb})
}
\\
\hline
24 & 52 & 18
\\
\hline
\end{tabular}
\label{tab:sesna_num_of_superficial_conv}
\end{center}
\end{table}

%%%%%%%%%%%%%%%%%%%%%%%%%%%%%%%%%%%%%%%%%%%%%%%%%%%%%%%%%%%%%%%%%%%%%%%%%

\Fig~\ref{fig:ConvLeftImpr1_woList} and
\Tab~\ref{tab:sesna_num_of_superficial_conv} show that
the stopping criterion (\ref{eqn:conv_judge2}\sq)
for the left-PCGS is inadequate.
Therefore,
the left-PCGS has a serious defect, in that superficial convergence can occur; we note that this
also occurs with other left-preconditioned algorithms \cite{itoh2010e}.

Just as information on \Fig~\ref{fig:ConvLeftImpr1_woList},
the numbers of problem cases
(not converged or not solved, breakdown, and superficial convergence)
are 172, 159, and 152
for the conventional PCGS, the left-PCGS, and the improved1 PCGS,
respectively.

\end{document}